\documentclass[a4paper]{article}

\usepackage[english]{babel}
\usepackage[T1]{fontenc}
\usepackage{amssymb}
\usepackage{amsthm}
\usepackage{makeidx}
\usepackage{color}
\usepackage{fullpage}
\usepackage{tikz}
\usepackage{mathtools}
\usepackage{multirow}
\usepackage[all]{xy}
\usepackage{comment}
\usepackage{authblk}
\usepackage[backend=biber,style=numeric,natbib=true]{biblatex}

\usepackage[shortlabels]{enumitem}
\usepackage{csquotes}
\usepackage{mathrsfs}

\usepackage[a4paper,top=3cm,bottom=2cm,left=3cm,right=3cm,marginparwidth=1.75cm]{geometry}

\usepackage{amsmath}
\usepackage{graphicx}
\usepackage[colorlinks=true, allcolors=blue]{hyperref}
\usepackage[all]{xy}

\providecommand{\MSC}[1]{\textbf{\textit{MSC---}} #1}

\newtheorem{defi}{Definition}[section]

\newtheorem{thm}[defi]{Theorem}
\newtheorem{coro}[defi]{Corollary}
\newtheorem{lem}[defi]{Lemma}

\newtheorem{pro}[defi]{Proposition}
\newtheorem{rmk}[defi]{Remark}
\newtheorem{prob}[defi]{Problem}

\title{The module structure of ideals in \\ degree $p$ extensions of $p$-adic fields}
\author[1]{Daniel Gil-Muñoz}
\affil[1]{Departament de Matem\`atiques, Universitat Polit\`ecnica de Catalunya, Campus Diagonal Besòs, Edifici A (EEBE), Av. Eduard Maristany 16, 08019 Barcelona, Spain \bigskip}

\date{}

\begin{document}
\maketitle

\begin{abstract}
Let $L/K$ be a degree $p$ extension of $p$-adic fields. We shall provide a complete characterization for the freeness of any fractional ideal $\mathfrak{p}_L^h$, $h\in\mathbb{Z}$, of $L$ as a module over its associated order $\mathfrak{A}_{L/K}(h)$ in the only Hopf-Galois structure on $L/K$, and an explicit description of $\mathfrak{A}_{L/K}(h)$. We show that there is no loss of generality in assuming that $0\leq h\leq p-1$ and that the normal closure of $L/K$ is totally ramified over $K$. If $L/K$ is maximally ramified, then $\mathfrak{A}_{L/K}(h)$ is the maximal $\mathcal{O}_K$-order in $H$ and $\mathfrak{p}_L^h$ is $\mathfrak{A}_{L/K}(h)$-free. The remaining cases correspond to typical extensions. For these, we provide an $\mathcal{O}_K$-basis of $\mathfrak{A}_{L/K}(h)$ and a characterization for the condition that $\mathfrak{p}_L^h$ is $\mathfrak{A}_{L/K}(h)$-free. In the non-exceptional cases, which had been partially studied in the Galois setting, the characterization depends on $h$ and the continued fraction expansion of $\frac{a}{p}$, where $a$ is the residue class mod $p$ of the (fractional) ramification jump $\ell$ of $L/K$. In the remaining exceptional case, where $L/K$ is almost maximally ramified and $a+1\leq h\leq p-1$, it consists in a numerical common minimizer condition.
\end{abstract}

\MSC{11S15,12F10,16T05}

\section{Introduction}

Let $p$ be an odd prime and let $L/K$ be an extension of $p$-adic fields. The fractional ideals of the ring of integers $\mathcal{O}_L$ in such an extension are of the form $\mathfrak{p}_L^h$, $h\in\mathbb{Z}$. If $L/K$ is Galois with group $G$, Galois module theory seeks to describe the structure of $\mathfrak{p}_L^h$ as a module over its associated order $\mathfrak{A}_{K[G]}(h)$ in $K[G]$, the Galois group algebra. This is defined as the maximal $\mathcal{O}_K$-order in $K[G]$ that acts on $\mathfrak{p}_L^h$. The most interesting case corresponds to $h=0$, where $\mathfrak{p}_L^h=\mathcal{O}_L$, the ring of integers of $L$, and it is studied as a module over its associated order in $K[G]$. The choice of the associated order as ground ring for the module structure is optimal in the sense that it is the only $\mathcal{O}_K$-order in $K[G]$ over which $\mathfrak{p}_L^h$ can possibly be free.

The picture provided by Galois module theory can be broadened with the use of Hopf-Galois theory, which is a generalization of Galois theory with the use of Hopf algebra. It starts with the notion of a Hopf-Galois structure on a finite field extension $L/K$: a pair formed by a $K$-Hopf algebra $H$ acting $K$-linearly on $L$ in the same fashion as a Galois group algebra does on a Galois extension. The extension $L/K$ is called Hopf-Galois if it admits some Hopf-Galois structure. Under this definition, every Galois extension is Hopf-Galois, as $K[G]$ together with its action on $L$ by endomorphisms is a Hopf-Galois structure, called the classical Galois structure. Also, for a Hopf-Galois extension $L/K$, one can consider the problem of characterizing the freeness of $\mathfrak{p}_L^h$ as a module over its associated order $\mathfrak{A}_H(h)$ in $H$, defined in an analogous way with respect to the Galois case. The study of this and adjacent questions belong to the so-called area of Hopf-Galois module theory.

In this paper, we provide a complete answer for this problem in the case in which $L/K$ is a degree $p$ extension of $p$-adic fields. Such an extension is always Hopf-Galois and admits a unique Hopf-Galois structure $H$. In the case that $L/K$ is Galois with group $G$, $H=K[G]$ is just its classical Galois structure. The freeness of the ring of integers in that case was completely characterized by F. Bertrandias, J.P. Bertrandias and M.J. Ferton \cite{bertrandiasbertrandiasferton1972,bertrandiasferton1972}. If $t$ is the ramification jump of such an extension, write $\frac{t}{p}=[a_0;a_1,\dots,a_n]$ for its continued fraction expansion.

\begin{thm}[F. Bertrandias, J.P. Bertrandias, M.J. Ferton]\label{thm:cyclicpfreeness} Let $L/K$ be a ramified cyclic degree $p$ extension of $p$-adic fields with ramification jump $t$, let $a$ be the remainder of $t$ mod $p$ and let $e$ be the ramification index of $K/\mathbb{Q}_p$.
\begin{enumerate}[label=\arabic*)]
    \item\label{thm:cyclicpfreeness1} If $a=0$, then $\mathcal{O}_L$ is $\mathfrak{A}_{L/K}$-free.
    \item\label{thm:cyclicpfreeness2} If $a\neq0$ and $t<\frac{pe}{p-1}-1$, then $\mathcal{O}_L$ is $\mathfrak{A}_{L/K}$-free if and only if $a\mid p-1$ (equivalently, $n\leq 2$).
    \item\label{thm:cyclicpfreeness3} If $a\neq0$ and $t\geq\frac{pe}{p-1}-1$, then $\mathcal{O}_L$ is $\mathfrak{A}_{L/K}$-free if and only if $n\leq 4$.
\end{enumerate}
\end{thm}

The statements \ref{thm:cyclicpfreeness1}-\ref{thm:cyclicpfreeness3} cover all possible cases because $t\leq\frac{pe}{p-1}$, with equality if and only if $a=0$. However, the techniques towards the proofs, especially for \ref{thm:cyclicpfreeness2} and \ref{thm:cyclicpfreeness3}, were rather opaque. Some subsequent works brought some light to this situation. Byott, Childs and Elder \cite{byottchildselder} introduced the notion of a scaffold as an object attached to an extension so as to derive criteria for freness in terms of numerical data. If $L/K$ is cyclic of degree $p$, they constructed a scaffold for $L/K$ in order to prove \ref{thm:cyclicpfreeness2} under the more restrictive condition that $t<\frac{pe}{p-1}-2$ (see Example 3.3 in that same reference). On the other hand, Del Corso, Ferri and Lombardo \cite{delcorsoferrilombardo} proved \ref{thm:cyclicpfreeness1} and \ref{thm:cyclicpfreeness2} using an alternative approach; namely, a study of the module index $[\mathcal{O}_L:\mathcal{O}_K[G]]$.

Still in the Galois case, the most complete study up to date for the structure of the ideals $\mathfrak{p}_L^h$ is due to Ferton \cite{fertonideals}. We restate her main result as follows.

\begin{thm}[\cite{fertonideals}, Th\'eor\`eme]\label{thm:fertonideals} Let $L/K$ be a ramified cyclic degree $p$ extension of $p$-adic fields. Let $t$ be the ramification jump of $L/K$ and let $a$ be the remainder of $t$ mod $p$.
\begin{enumerate}[label=\arabic*)]
    \item If $a=0$, $\mathfrak{p}_L^h$ is $\mathfrak{A}_{L/K}(h)$-free for all $0\leq h\leq p-1$.
    \item If $a=1$, $h\in\{0,1\}$ and $t<\frac{pe}{p-1}-2$, then $\mathfrak{p}_L^h$ is $\mathfrak{A}_{L/K}(h)$-free.
    \item If $a\neq 1$, $0\leq h\leq a$ and $t<\frac{pe}{p-1}-1$, then $\mathfrak{p}_L^h$ is $\mathfrak{A}_{L/K}(h)$-free if and only if:
    \begin{enumerate}[label=\alph*)]
        \item\label{thm:fertonideals3a} $h\in\{a-a_n,a\}$, if $n$ is even.
        \item\label{thm:fertonideals3b} $a-\frac{a_n}{2}\leq h\leq a$, if $n$ is odd.
    \end{enumerate}
    \item If $a+1\leq h\leq p-1$ and $t<\frac{pe}{p-1}-2$, then $\mathfrak{p}_L^h$ is $\mathfrak{A}_{L/K}(h)$-free if and only if $a=1$ and $h>\frac{p+1}{2}$.
\end{enumerate}
\end{thm}

Note that this result does not cover any ideal $\mathfrak{p}_L^h$ when $\frac{pe}{p-1}-1\leq t<\frac{pe}{p-1}$, and only covers those with $0\leq h\leq a$ and $a\neq1$ when $\frac{pe}{p-1}-2\leq t<\frac{pe}{p-1}-1$.

On the Hopf-Galois side, an explicit description of the Hopf-Galois structure $H$ on $L/K$ is known. It was first provided by Childs \cite{childs1989}, who used the description of rank $p$ Hopf algebras by Tate and Oort \cite{tateoort} so as to find a $K$-algebra generator $\Psi$ of $H$. This description is replayed at \cite{elder2018} for the so called typical degree $p$ extensions (in which case the ring of integers is not generated by a $p$-th root of the uniformizer) and proved as valid also for the prime characteristic case. In addition, the author proved that $\Psi^p=\lambda\Psi$ for a specific $\lambda\in K$ (see \cite[Corollary 5.9]{gildegpgeneral}).

As for freeness, the author extended Theorem \ref{thm:cyclicpfreeness} to an arbitrary degree $p$ extension of $p$-adic fields $L/K$, first to the case in which the normal closure is dihedral \cite{gildegpdihedral}, and afterwards in the general case \cite{gildegpgeneral}. In this characterization, the role of $t$ is played by the ramification jump $\ell$ of $L/K$. If $L/K$ is not Galois, $\ell$ is defined through the Herbrand function of $L/K$ (see for instance \cite[\textsection 2.2.1]{gildegpgeneral}). In that case, $\ell$ is not a rational integer but a $p$-adic integer. Thus, $a$ can be defined as the residue class of $\ell$ mod $p$, and $\ell\leq\frac{pe}{p-1}$, with equality if and only if $a=0$. Write $$\frac{a}{p}=[0;a_1,\dots,a_n]$$ for the continued fraction expansion of $\frac{a}{p}$. For the generalization of Theorem \ref{thm:cyclicpfreeness} \ref{thm:cyclicpfreeness3}, the condition characterizing freeness is $n\leq 4$. On the other hand, Elder used the theory of scaffolds to give a characterization for the $\mathfrak{A}_{L/K}(h)$-freeness of $\mathfrak{p}_L^h$ under some conditions (see \cite[Theorem 3.5]{elder2018}). Namely, there is always freeness for $h\equiv a\,(\mathrm{mod}\,p)$, and if $\ell<\frac{pe}{p-1}-2$, he provided characterizations for the cases $h\equiv0\,(\mathrm{mod}\,p)$ and $h\equiv a+1\,(\mathrm{mod}\,p)$.

We shall provide a complete characterization for the $\mathfrak{A}_{L/K}(h)$-freeness of each fractional ideal $\mathfrak{p}_L^h$ and an explicit $\mathcal{O}_K$-basis of $\mathfrak{A}_{L/K}(h)$, covering all gaps in previous works. There are two significant simplifications to this problem. One of them is that the $\mathfrak{A}_{L/K}(h)$-freeness of $\mathfrak{p}_L^h$ is periodic on $h$ (with period $p$), so one can assume without loss of generality that $0\leq h\leq p-1$. On the other hand, we can also assume without loss of generality that the normal closure $\widetilde{L}$ of $L/K$ is totally ramified over $K$; otherwise we would replace $L/K$ by the extension $L'/K'$ of the inertia fields of $\widetilde{L}/L$ and $\widetilde{L}/K$. Our first main result is as follows.


\begin{thm}\label{thm:freenessideals} Let $L/K$ be a degree $p$ extension of $p$-adic fields with normal closure $\widetilde{L}$, and assume that $\widetilde{L}/K$ is totally ramified. Call $e=e(K/\mathbb{Q}_p)$. Let $\ell$ be the ramification break of $L/K$ and let $a$ be its residue class modulo $p$. Fix $0\leq h\leq p-1$. Then, the following statements hold.
\begin{enumerate}[label=\arabic*)]
    \item\label{thm:freenessideals1} If $a=0$, then $\mathfrak{A}_{L/K}(h)$ is the maximal $\mathcal{O}_K$-order in $H$ and $\mathfrak{p}_L^h$ is $\mathfrak{A}_{L/K}(h)$-free.
    \item\label{thm:freenessideals2} If $1\leq h\leq a$, then $\mathfrak{p}_L^h$ is $\mathfrak{A}_{L/K}(h)$-free if and only if:
    \begin{enumerate}[label=\alph*)]
        \item\label{thm:freenessideals2a} $h\in\{a-a_n,a\}$, if $n$ is even.
        \item\label{thm:freenessideals2b} $a-\frac{a_n}{2}\leq h\leq a$, if $n$ is odd.
    \end{enumerate}
    \item\label{thm:freenessideals3} If $a+1\leq h\leq p-1$ and $\ell<\frac{pe}{p-1}-1$, then $\mathfrak{p}_L^h$ is $\mathfrak{A}_{L/K}(h)$-free if and only if $a=1$ and $h>\frac{p+1}{2}$.
\end{enumerate}
\end{thm}

The case $a\neq0$ and $h=0$ is not covered by Theorem \ref{thm:freenessideals}. However, that situation is precisely when $\mathfrak{p}_L^h=\mathcal{O}_L$, and \cite[Theorem 1.2 4)]{gildegpgeneral} corresponds to a characterization of its $\mathfrak{A}_H$-freeness. If one takes $h=0$ in the statements \ref{thm:freenessideals2a}-\ref{thm:freenessideals2b}, then one recovers the condition that $a\mid p-1$, and this condition is sufficient but not necessary. This means that the continued fraction phenomena arising at the behavior of the $\mathfrak{A}_{L/K}$-freeness of $\mathcal{O}_L$ in Theorem \ref{thm:cyclicpfreeness} \ref{thm:cyclicpfreeness3} (and its Hopf-Galois generalization) does not occur for $1\leq h\leq a$.

If we specialize Theorem \ref{thm:freenessideals} to the case in which $L/K$ is Galois, we can complete some of the gaps left by Theorem \ref{thm:fertonideals}, as follows:
\begin{itemize}
    \item If $a=1$, $h\in\{0,1\}$ and $t\geq\frac{pe}{p-1}-2$, then $\mathfrak{p}_L^h$ is $\mathfrak{A}_{L/K}(h)$-free.
    \item If $a\neq 1$, $0\leq h\leq a$ and $t\geq\frac{pe}{p-1}-1$, then $\mathfrak{p}_L^h$ is $\mathfrak{A}_{L/K}(h)$-free if and only if the same two conditions \ref{thm:fertonideals3a}-\ref{thm:fertonideals3b} hold.
    \item If $a+1\leq h\leq p-1$ and $\frac{pe}{p-1}-2\leq t<\frac{pe}{p-1}-1$, then $\mathfrak{p}_L^h$ is $\mathfrak{A}_{L/K}(h)$-free if and only if $a=1$ and $h>\frac{p+1}{2}$.
\end{itemize}

The other uncovered case by Theorem \ref{thm:freenessideals} is $a+1\leq h\leq p-1$ and $\ell\geq\frac{pe}{p-1}-1$, which had not been studied in the Galois case neither. We shall state a characterization for freeness also for this case, for which we will need a few preparations. Let $c$ the remainder of $r\ell$ mod $r$ and $b=\ell-\frac{pc}{r}$. For each $0\leq k\leq p-1$, call $$d(k;h)=\Big\lfloor\frac{a'+kb-h}{p}\Big\rfloor,$$ where $a'=a$ if $0\leq h\leq a$ and $a'=a+p$ if $a+1\leq h\leq p-1$. The key parameters governing freeness are $$W(i,j;h)=\begin{cases}
    d(i+j;h)-d(j;h) & \hbox{if }i+j\leq p-1, \\
    v_K(\lambda)+d(i+j+1-p;h)-d(j;h) & \hbox{if }i+j\geq p,
\end{cases}$$ where $0\leq i,j\leq p-1$. Call $$\widetilde{w}(i;h)=\mathrm{min}\{W(i,j;h)\,\mid\,0\leq j\leq p-1\}.$$ 

\begin{thm}\label{thm:freenessexceptionalcase} Let $L/K$ be a typical degree $p$ extension of $p$-adic fields with normal closure $\widetilde{L}$, and suppose that $\widetilde{L}/K$ is totally ramified. If $a+1\leq h\leq p-1$ and $\ell\geq\frac{pe}{p-1}-1$, then $\mathfrak{p}_L^h$ is $\mathfrak{A}_{L/K}(h)$-free if and only if there is $0\leq s\leq p-1$ such that $$W(i,s;h)=\widetilde{w}(i;h)$$ for all $0\leq i\leq p-1$.
\end{thm}

Consider now the problem of describing $\mathfrak{A}_{L/K}(h)$. In the case that $a=0$, it is known that associated order $\mathfrak{A}_{L/K}$ in the ring of integers is the maximal $\mathcal{O}_K$-order in $H$. We shall see that $\mathfrak{A}_{L/K}(h)=\mathfrak{A}_{L/K}$, from which Theorem \ref{thm:freenessideals} \ref{thm:freenessideals1} will follow. In the typical case, the elements $\widetilde{w}(i;h)$ can be also used so as to determine an $\mathcal{O}_K$-basis of $\mathfrak{A}_{L/K}(h)$.

\begin{thm}\label{thm:assocordtypical} Let $L/K$ be a typical degree $p$ extension of $p$-adic fields with totally ramified normal closure. Then, $\mathfrak{A}_{L/K}(h)$ has $\mathcal{O}_K$-basis $$\{\pi_K^{-\widetilde{w}(i;h)}\Psi^i\}_{i=0}^{p-1}.$$
\end{thm}

If $0\leq h\leq a$ or $\ell<\frac{pe}{p-1}-1$, then $\mathfrak{A}_{L/K}(h)$ has $\mathcal{O}_K$-basis $\{\pi_K^{-w(i;h)}\Psi^i\}_{i=0}^{p-1}$, where $$w(i;h)=\mathrm{min}\{d(i+j;h)-d(j;h)\,\mid\,0\leq j\leq p-1-i\}$$ (see Corollary \ref{coro:criteriamidtransl}). This is coherent with Theorem \ref{thm:assocordtypical}, as we shall see that $\widetilde{w}(i;h)=w(i;h)$ under these restrictions. 

\subsection*{Structure of the paper}

This paper is organized as follows. In Section \ref{sec:hgmfractid}, we shall view some standard facts on the associated order $\mathfrak{A}_H(h)$ in the case of an arbitrary extension $L/K$ of $p$-adic fields, and we will establish various characterizations for the $\mathfrak{A}_H(h)$-freeness of $\mathfrak{p}_L^h$ that will be used along the paper. It is in this part where the periodicity of $h$ for freeness is proved. In Section \ref{sec:degpsetup}, we will specialize to the case in which $L/K$ has degree $p$, and we shall establish the setup that will be needed. In addition, we will go over the already known description of $H$ as a $K$-algebra generated by $\Psi$. In Section \ref{sec:redtotram}, the case in which $\widetilde{L}/K$ is not totally ramified is treated, reducing it to the totally ramified case. On the other hand, the case $a=0$ is treated in Section \ref{sec:maxram}, where Theorem \ref{thm:freenessideals} \ref{thm:freenessideals1} is proved. 

As already mentioned, the case that $a\neq0$ corresponds to typical extensions. In Section \ref{sec:typicaldegp}, we recall Elder's work on them; namely, a characterization in terms of a defining equation, the construction of a scaffold on $L/K$ and its application so as to derive sufficient conditions for freeness. We shall see that Elder's conclusion can actually be reached under the assumption that $\ell<\frac{pe}{p-1}-1$ if $0\leq h\leq a$. It is in Section \ref{sec:equaldiwi} where we assume that $h>0$ and prove that the sufficient conditions provided by the scaffold are also necessary, discarding the continued fraction phenomena for this case. In Section \ref{sec:equalassocorders}, we shall rewrite these conditions in the form in which they appear in Theorem \ref{thm:freenessideals}, by characterizing when $\mathfrak{A}_{L/K}(h)$ coincides with an $\mathcal{O}_K$-order $\mathfrak{A}_{\pi_L^{a'}}(h)$. The statements \ref{thm:freenessideals2} and \ref{thm:freenessideals3} will follow from the results in these two sections. The exceptional case $a+1\leq h\leq p-1$ and $\ell\geq\frac{pe}{p-1}-1$ is treated in Section \ref{sec:exceptionalcase}, where we prove Theorems \ref{thm:freenessexceptionalcase} and \ref{thm:assocordtypical}.

\section{Preliminaries}

Unless specified otherwise, all our field extensions will be assumed to be finite.

\subsection{Hopf-Galois extensions}

Let $L/K$ be a finite field extension. A Hopf-Galois structure on $L/K$ is a pair formed by a $K$-Hopf algebra $H$ and a $K$-linear map $\cdot\colon H\otimes_KL\longrightarrow L$ such that $L$ is an $H$-module algebra with respect to $\cdot$ (see \cite[Chapter 2, Section 2.1, Page 39]{underwood2015}) and the map $j\colon L\otimes_KH\longrightarrow\mathrm{End}_K(L)$ defined by $$j(x\otimes h)(y)=x(h\cdot y),\quad x,y\in L,\,h\in H$$ is a $K$-linear isomorphism. If $L/K$ admits some Hopf-Galois structure, we say that $L/K$ is Hopf-Galois. We will also say that $L/K$ is $H$-Galois. We shall use the label $H$ for both the Hopf-Galois structure and for its underlying Hopf algebra. It is obvious that $H$ endows $L$ with module structure. It holds that $L$ is $H$-free of rank one (see \cite[(2.16)]{childs}). This is a generalization of Normal Basis Theorem to Hopf-Galois extensions. Consequently, any element $\theta\in H$ with $L=H\cdot\theta$ will be called an $H$-normal basis generator for $L$.

Suppose that $L/K$ is finite and separable. Let $\widetilde{L}$ be the normal closure of $L/K$, $G=\mathrm{Gal}(\widetilde{L}/K)$ and $G'=\mathrm{Gal}(\widetilde{L}/L)$. Let $\lambda\colon G\longrightarrow\mathrm{Perm}(G/G')$ be the map defined by $\lambda(g)(\overline{h})=\overline{gh}$, which is usually known as the left translation map for $G$. It is known to be a group monomorphism (see \cite[(6.6)]{childs}). 

Greither and Pareigis \cite{greitherpareigis} proved that the set of Hopf-Galois structures on $L/K$ is in bijective correspondence with the set of regular subgroups of $\mathrm{Perm}(G/G')$ that are normalized by $\lambda(G)$. Here a subgroup of $\mathrm{Perm}(G/G')$ is said to be regular if its action on $G/G'$ is simply transitive. 

Moreover, if $N$ is a regular subgroup normalized by $\lambda(G)$, the corresponding Hopf-Galois structure on $L/K$ is given by the $K$-Hopf algebra $$\widetilde{L}[N]^G=\{h\in\widetilde{L}[N]\,\mid\,g*h=h\hbox{ for all }g\in G\},$$ where $$g*\Big(\sum_{\eta\in N}h_{\eta}\eta\Big)=\sum_{\eta\in N}g(h_{\eta})\lambda(g)\eta\lambda(g)^{-1},\quad g\in G,\,h_{\eta}\in\widetilde{L},$$ and the action of $h=\sum_{\eta\in N}h_{\eta}\eta\in\widetilde{L}[N]^G$ on $L$ given by $$h\cdot x=\sum_{\eta\in N}h_{\eta}\eta^{-1}(\overline{1_G})(x),\quad x\in L.$$ The expression $\eta^{-1}(\overline{1_G})(x)$ denotes a representative of $\eta^{-1}(\overline{1_G})\in G/G'$ evaluated on $x$. Since $L$ is fixed by the elements of $G'$, the definition is correct.

We say that the extension $L/K$ is almost classically Galois if $G'$ admits some normal complement $J$ in $G$, so that $G=J\rtimes G'$. In that case, $\lambda(J)$ is a regular subgroup normalized by $\lambda(G)$, so it gives rise to a Hopf-Galois structure on $L/K$. In this case, the description of its underlying Hopf algebra is slightly simplified to $$H=M[J]^{G'},$$ where $M\coloneqq\widetilde{L}^J$ (see for instance \cite[Proof of Theorem 3.1]{KKTU19}).

An integer number $n$ is said to be Burnside if it is coprime with $\varphi(n)$, its image by the Euler totient function $\varphi$. Byott's uniqueness theorem states that if $L/K$ is a separable extension with Burnside degree, then it admits a unique Hopf-Galois structure $H$, and in addition $L/K$ is almost classically Galois. 

A particularly simple case for the search of Hopf-Galois structures corresponds to extensions $L/K$ with prime degree. Acquiring Greither-Pareigis' notation, such an extension is Hopf-Galois if and only if $G$ is solvable (see \cite[(7.5)]{childs}). Since any prime number is Burnside, Byott's uniqueness theorem yields that prime degree Hopf-Galois separable extensions admit a unique Hopf-Galois structure.

\subsection{Local fields with mixed characteristic and their ramification}

Let $p$ be an odd prime. By a $p$-adic field we mean a finite extension of the field $\mathbb{Q}_p$ of $p$-adic numbers. For a $p$-adic field $F$, we will use the following notations:
\begin{itemize}
    \item $\mathcal{O}_F$ is the ring of integers or valuation ring of $F$.
    \item $\mathfrak{p}_F$ is the prime ideal of $\mathcal{O}_F$.
    \item $v_F$ is the valuation of $F$.
    \item $\pi_F$ is a uniformizer of $F$.
    \item $\overline{F}=\mathcal{O}_F/\mathfrak{p}_F$ is the residue field of $F$.
\end{itemize}

We know that $\overline{F}$ is a field with characteristic $p$. For this reason, $p$-adic fields are commonly known as local fields with mixed characteristic. The ring $\mathcal{O}_F$ is a discrete valuation ring, as its only prime ideal is $\mathfrak{p}_F$. For this reason, its ideals are of the form $\mathfrak{p}_F^h$ with $h\in\mathbb{Z}_{\geq0}$. We follow the convention that $\mathfrak{p}_F^0=\mathcal{O}_F$. In addition, if $h<0$, we shall write $\mathfrak{p}_F^h$ for the elements of $F$ with valuation at least $h$. This is a fractional ideal of $\mathcal{O}_F$. In fact, the fractional ideals of $\mathcal{O}_F$ are $\mathfrak{p}_F^h$, where $h\in\mathbb{Z}$.

Let $E/K$ be an extension of $p$-adic fields. Then, $\mathfrak{p}_K\mathcal{O}_E$ is an ideal of $\mathcal{O}_E$, so there is some $e\in\mathbb{Z}_{\geq1}$ such that $\mathfrak{p}_K\mathcal{O}_E=\mathfrak{p}_E^e$. The integer $e$ is by definition the ramification index of $E/K$, denoted by $e(E/K)$. Also, it is known that there is a field monomorphism $\overline{K}\longrightarrow\overline{E}$. The residue class degree is defined as the degree $f(E/K)\coloneqq[\overline{E}:\overline{K}]$. It is known that $[E:K]=e(E/K)f(E/K)$ (see for instance \cite[Chapter 24, Theorem 1]{lorenz}). We say that $E/K$ is unramified if $e(E/K)=1$, and ramified otherwise. If $e(E/K)=[E:K]$, we will say that $E/K$ is totally ramified.

Next, suppose that $E/K$ is Galois and let $G=\mathrm{Gal}(E/K)$. Since $\overline{E}/\overline{K}$ is separable, there is some $x\in\mathcal{O}_E$ such that $\mathcal{O}_E=\mathcal{O}_K[x]$ (see \cite[Chapter III, \textsection6, Proposition 12]{serre}). For an integer number $i\geq-1$, the $i$-th ramification group of $E/F$ is $G_{-1}=G$ if $i=-1$ and $$G_i=\{g\in G\,\mid\,v_E(g(x)-x)\geq i+1\}$$ if $i\geq0$. It is clear that $G_i\supseteq G_{i+1}$ for every $i\geq-1$ and that $G_i$ is trivial for $i$ large enough (namely, $i\geq\mathrm{max}_{g\in G}(g(x)-x)$). Therefore, $\{G_i\}_{i=-1}^{\infty}$ is a filtration of $G$, called the chain of ramification groups of $E/K$. A \textbf{ramification jump} for $E/K$ is an integer $t\geq-1$ such that $G_t\neq G_{t+1}$. This concept is also known as ramification break or ramification number in literature. Since $\mid G_0\mid=e(E/K)$, $E/K$ is unramified if and only if $G_0=\{1\}$, and it is totally ramified if and only if $G_0=G$. It is also known that $G_0$ is solvable \cite[Chapter IV, \textsection2, Corollary 5]{serre} and $G_1$ is the $p$-Sylow subgroup of $G$ \cite[Chapter IV, \textsection2, Corollary 3]{serre}.

The Herbrand function of $E/K$ is the function $\varphi_{E/K}\colon\mathbb{R}_{\geq-1}\longrightarrow\mathbb{R}_{\geq-1}$ defined as $$\varphi_{E/K}(u)=\int_0^u\frac{1}{[G_0:G_t]}dt,$$ where $G_u=G_{\lfloor u\rfloor}$ for all $u\in\mathbb{R}$ (see \cite[Chapter IV, \textsection3]{serre}). It holds that $\varphi_{E/K}$ is continuous and strictly increasing (see \cite[Chapter IV, \textsection3, Proposition 12 a)]{serre}), therefore bijective; call $\psi_{E/K}$ its inverse. If $N$ is an intermediate field of $E/K$, by \cite[Chapter IV, \textsection3, Proposition 15]{serre}, we have transitivity formulas \begin{equation}\label{eq:transitherbrand}
    \varphi_{E/K}=\varphi_{N/K}\varphi_{E/N},\quad \psi_{E/K}=\psi_{E/N}\psi_{N/K}.
\end{equation} 

Let $L/K$ be a (non-necessarily Galois) extension of $p$-adic fields with normal closure $\widetilde{L}\subseteq E$. We can extend the definition of Herbrand function from the Galois case to the extension $L/K$ as follows: $$\varphi_{L/K}\coloneqq\varphi_{E/K}\psi_{E/L}.$$ This definition does not depend on the choice of Galois extension $E/K$ containing $L$ because of the relations \eqref{eq:transitherbrand}. The image of a ramification jump for $E/K$ by $\varphi_{L/K}$ will be referred to as a ramification jump for $L/K$. This might not be an integer number, but in any case it is a $p$-adic integer.

\section{Hopf-Galois module structure of fractional ideals}\label{sec:hgmfractid}

Let $L/K$ be an $H$-Galois extension of $p$-adic fields. The problem that is usually studied in Hopf-Galois module theory is whether $\mathcal{O}_L$ is free as a module over its associated order $\mathfrak{A}_H$ in $H$, which is defined as the maximal $\mathcal{O}_K$-order in $H$ that acts on $\mathcal{O}_L$. We shall work with the following more general object.

\begin{defi}\label{def:assocordideal} Let $h\in\mathbb{Z}$. The associated order of the fractional ideal $\mathfrak{p}_L^h$ in $H$ is defined as $$\mathfrak{A}_H(h)\coloneqq\{\lambda\in H\,\mid\,\lambda\cdot\mathfrak{p}_L^h\subset\mathfrak{p}_L^h\}.$$
\end{defi}

For $h=0$, we have that $\mathfrak{A}_H(0)=\mathfrak{A}_H$, the usual associated order of $\mathcal{O}_L$ in $H$. In analogy to that case, we can consider the problem of studying the $\mathfrak{A}_H(h)$-module structure of $\mathfrak{p}_L^h$. Even though the properties of $\mathfrak{A}_H(h)$ that provide sense to such an study are well known, we prove them explicitly here for expository purposes.

Recall that $L/K$ being $H$-Galois is defined as the bijectivity of the canonical map $j\colon L\otimes_KH\longrightarrow\mathrm{End}_K(L)$. Its restriction to $H$ is the $K$-linear representation $\rho\colon H\longrightarrow\mathrm{End}_K(L)$ defined by $\rho(h)(x)=h\cdot x$ for $h\in H$ and $x\in L$. Since $j$ is bijective, $\rho$ is injective. The map $\rho$ identifies the action of $H$ in terms of endomorphisms, and the associated order $\mathfrak{A}_H(h)$ has a precise meaning in these terms.

\begin{lem}\label{lem:assocordend} Let $h\in\mathbb{Z}$. Then $$\rho(\mathfrak{A}_H(h))=\rho(H)\cap\mathrm{End}_{\mathcal{O}_K}(\mathfrak{p}_L^h).$$
\end{lem}
\begin{proof}
Given $\gamma\in H$, we have $\gamma\in\mathfrak{A}_H(h)$ if and only if $\rho(\gamma)(\mathfrak{p}_L^h)=\gamma\cdot\mathfrak{p}_L^h\subset\mathfrak{p}_L^h$. This is the same as $\rho(\gamma)\in\mathrm{End}_{\mathcal{O}_K}(\mathfrak{p}_L^h)$, and of course $\rho(\gamma)\in\rho(H)$. Then, the equality follows.
\end{proof}

With the help of this lemma, we can prove that the associated order $\mathfrak{A}_H(h)$ is actually an $\mathcal{O}_K$-order.

\begin{pro}\label{pro:assocordproperties} Fix $h\in\mathbb{Z}$ and consider the fractional ideal $\mathfrak{p}_L^h$ of $\mathcal{O}_L$.
\begin{enumerate}[label=\arabic*)]
    \item\label{pro:assocordproperties1} $\mathfrak{A}_H(h)$ is an $\mathcal{O}_K$-order in $H$.
    \item\label{pro:assocordproperties2} If $\mathfrak{A}$ is an $\mathcal{O}_K$-order in $H$ over which $\mathfrak{p}_L^h$ is $\mathfrak{A}$-free, then $\mathfrak{A}=\mathfrak{A}_H(h)$.
\end{enumerate}
\end{pro}
\begin{proof}
\begin{enumerate} 
    \item Note that $\mathrm{End}_{\mathcal{O}_K}(\mathfrak{p}_L^h)$ is naturally identified with the space of matrices $\mathcal{M}_n(\mathcal{O}_K)$. Then, $\mathrm{End}_{\mathcal{O}_K}(\mathfrak{p}_L^h)$ is $\mathcal{O}_K$-finitely generated. Since in addition $\mathcal{O}_K$ is noetherian and $\rho(\mathfrak{A}_H(h))\subset\mathrm{End}_{\mathcal{O}_K}(\mathfrak{p}_L^h)$ by Lemma \ref{lem:assocordend}, we deduce that $\rho(\mathfrak{A}_H(h))$ is $\mathcal{O}_K$-finitely generated. Now, as $\rho$ injective, it maps isomorphically $\mathfrak{A}_H(h)$ over its image, so $\mathfrak{A}_H(h)$ is $\mathcal{O}_K$-finitely generated.
    Let $\{\gamma_i\}_{i=1}^n$ be a $K$-basis of $H$. Given $1\leq i\leq n$, $\gamma_i\cdot\mathfrak{p}_L^h\subset L$. Since $\mathfrak{p}_L^h$ is a finitely generated torsion-free $\mathcal{O}_K$-module, there is $m_i\in\mathbb{Z}$ such that $\pi_K^{m_i}\gamma_i\cdot\mathfrak{p}_L^h\subset\mathfrak{p}_L^h$, so $\pi_K^{m_i}\gamma_i\in\mathfrak{A}_H(h)$. Thus, $\{\pi_K^{m_i}\gamma_i\}_{i=1}^n$ is a $K$-basis of $H$ contained in $\mathfrak{A}_H(h)$.
    \item Suppose that $\mathfrak{p}_L^h$ is $\mathfrak{A}$-free, say $\mathfrak{p}_L^h=\mathfrak{A}\cdot\theta$ with $\theta\in\mathfrak{p}_L^h$. By definition of associated order, $\mathfrak{A}\subseteq\mathfrak{A}_H(h)$, so we just need to prove the reverse inclusion. Let $\gamma\in\mathfrak{A}_H(h)$. Since $\gamma\cdot\theta\in\mathfrak{p}_L^h$, there is some $\lambda\in\mathfrak{A}$ such that $\gamma\cdot\theta=\lambda\cdot\theta$. Since $L$ is $H$-free, we have that $\lambda=\gamma$. We see that $\gamma\in\mathfrak{A}$, proving the desired equality.
\end{enumerate}
\end{proof}

Hence, as in the case of $\mathfrak{A}_H$, $\mathfrak{A}_H(h)$ is the only $\mathcal{O}_K$-order in $H$ over which $\mathfrak{p}_L^h$ is possibly free. There is a significant simplification that can be done to this problem.

\begin{lem}\label{lem:equivfreenesses} Call $e\coloneqq e(L/K)$ and let $h_1,h_2\in\mathbb{Z}$ with $h_1\equiv h_2\,(\mathrm{mod}\,e)$. Then $\mathfrak{A}_{H}(h_1)=\mathfrak{A}_H(h_2)$ and $\mathfrak{p}_L^{h_1}$ is $\mathfrak{A}_H(h_1)$-free if and only if $\mathfrak{p}_L^{h_2}$ is $\mathfrak{A}_H(h_2)$-free.
\end{lem}
\begin{proof}
Let $k\coloneqq\frac{h_1-h_2}{e}\in\mathbb{Z}$, so that $h_1=ke+h_2$. We see that $\mathfrak{p}_L^{h_1}=\pi_K^k\mathfrak{p}_L^{h_2}$. Then, the map $\psi\colon\mathfrak{p}_L^{h_2}\longrightarrow\mathfrak{p}_L^{h_1}$ defined by $\psi(x)=\pi_K^kx$ is an $\mathcal{O}_K$-linear isomorphism. In addition, given $\gamma\in H$ and $x\in\mathfrak{p}_L^{h_2}$, $$\gamma\cdot\psi(x)=\gamma\cdot(\pi_K^kx)=\pi_K^k\gamma\cdot x=\psi(\gamma\cdot x).$$ In other words, $\psi$ is equivariant with respect to the Hopf action.

Now, given $\gamma\in H$, \begin{equation*}
    \begin{split}
        \gamma\in\mathfrak{A}_H(h_1)&\Longleftrightarrow\gamma\cdot\mathfrak{p}_L^{h_1}\subseteq\mathfrak{p}_L^{h_1}\\&\Longleftrightarrow\gamma\cdot\psi(\mathfrak{p}_L(h_2))\subset\psi(\mathfrak{p}_L^{h_2})
        \\&\Longleftrightarrow\psi(\gamma\cdot\mathfrak{p}_L^{h_2})\subset\psi(\mathfrak{p}_L^{h_2})
        \\&\Longleftrightarrow\gamma\cdot\mathfrak{p}_L^{h_2}\subset\mathfrak{p}_L^{h_2}
        \\&\Longleftrightarrow\gamma\in\mathfrak{A}_H(h_2).
    \end{split}
\end{equation*} We obtain that $\mathfrak{A}_H(h_1)=\mathfrak{A}_H(h_2)$.

On the other hand, for $\beta\in\mathfrak{p}_L^{h_2}$, $\mathfrak{p}_L^{h_2}=\mathfrak{A}_H(h_2)\cdot\beta$ if and only if $\mathfrak{p}_L^{h_1}=\mathfrak{A}_H(h_1)\cdot\psi(\beta)$, because $\psi$ is equivariant and maps $\mathfrak{p}_L^{h_2}$ isomorphically to $\mathfrak{p}_L^{h_1}$. Hence, the last sentence of the statement follows.
\end{proof}

Thus, in order to study the $\mathfrak{A}_H(h)$ freeness of $\mathfrak{p}_L^h$ for all $h\in\mathbb{Z}$, it is equivalent to do so just for $0\leq h\leq e(L/K)-1$. Namely:

\begin{prob} Let $L/K$ be an $H$-Galois extension of $p$-adic fields and let $h\in\mathbb{Z}$ with $0\leq h\leq e(L/K)-1$. Find a necessary and sufficient condition for $\mathfrak{p}_L^h$ being $\mathfrak{A}_H(h)$-free.
\end{prob}

This problem has been considered in \cite{byottchildselder} for totally ramified degree $p^n$ extensions of local fields, where $e(L/K)=p^n$, by means of the theory of scaffolds, which we recall in Section \ref{sec:scaffoldp} for $n=1$. 

For $\theta\in\mathfrak{p}_L^h$, denote $$\mathfrak{A}_{\theta}(h)=\{\lambda\in H\,\mid\,\lambda\cdot\theta\in\mathfrak{p}_L^h\}.$$ It is immediate that $\mathfrak{A}_H(h)\subset\mathfrak{A}_{\theta}(h)$. 

\begin{lem}\label{lem:atheta} Let $\theta\in\mathfrak{p}_L^h$ be an $H$-normal basis generator for $L$. Then:
\begin{enumerate}[label=\arabic*)]
    \item\label{lem:atheta1} $\mathfrak{p}_L^h=\mathfrak{A}_{\theta}(h)\cdot\theta$.
    \item\label{lem:atheta2} $\mathfrak{A}_{\theta}(h)$ is an $\mathcal{O}_K$-lattice in $H$.
    \item\label{lem:atheta3} $\mathfrak{A}_{\theta}(h)$ is an $\mathfrak{A}_H(h)$-fractional ideal.
\end{enumerate}
\end{lem}
\begin{proof}
\begin{enumerate}[label=\arabic*)]
    \item We know from the definition that $\mathfrak{A}_{\theta}(h)\cdot\theta\subseteq\mathfrak{p}_L^h$. Let $x\in\mathfrak{p}_L^h\subset L$. Since $\theta$ is an $H$-normal basis generator, there is $\gamma\in H$ such that $x=\gamma\cdot\theta$. Since $\gamma$ acting on $\theta$ falls in $\mathfrak{p}_L^h$, we deduce that $\gamma\in\mathfrak{A}_{\theta}(h)$, proving that $\mathfrak{p}_L^h=\mathfrak{A}_{\theta}(h)\cdot\theta$.
    \item The assumption that $\theta$ is an $H$-normal basis generator for $L$ yields that the map $\Phi_{\theta}\colon H\longrightarrow L$ defined by $\Phi_{\theta}(\gamma)=\gamma\cdot\theta$ is a $K$-linear isomorphism. Now, the equality at \ref{lem:atheta1} is equivalent to $\Phi_{\theta}(\mathfrak{A}_{\theta}(h))=\mathfrak{p}_L^h$. Since $\mathfrak{p}_L^h$ is an $\mathcal{O}_K$-lattice in $L$, \ref{lem:atheta2} follows.
    \item Given $\gamma\in\mathfrak{A}_H(h)$ and $\lambda\in\mathfrak{A}_{\theta}(h)$, $$(\gamma\lambda)\cdot\theta=\gamma\cdot(\lambda\cdot\theta)\in\gamma\cdot\mathfrak{p}_L^h\subset\mathfrak{p}_L^h,$$ so $\gamma\lambda\in\mathfrak{A}_{\theta}(h)$, proving that $\mathfrak{A}_{\theta}(h)$ is an $\mathfrak{A}_H(h)$-module. Since both $\mathfrak{A}_H(h)$ and $\mathfrak{A}_{\theta}(h)$ are full $\mathcal{O}_K$-lattices in $H$, there is $m\in\mathbb{Z}_{\geq0}$ such that $\pi_K^m\mathfrak{A}_{\theta}(h)\subseteq\mathfrak{A}_H(h)$. As $\pi_K^m$ is not a zero divisor of $\mathcal{O}_K$, this finishes the proof of \ref{lem:atheta3}.
\end{enumerate}
\end{proof}

Recall that, for fractional ideals $I$ and $J$ of a ring $R$, the ideal quotient $(I:J)$ is defined as $$(I:J)\coloneqq\{x\in R\,\mid\,xI\subseteq J\}.$$ We can characterize $\mathfrak{A}_H(h)$ by means of this object.

\begin{coro}\label{coro:quotideal} Let $\theta\in\mathfrak{p}_L^h$ be an $H$-normal basis generator for $L$. Then $$\mathfrak{A}_H(h)=(\mathfrak{A}_{\theta}(h):\mathfrak{A}_{\theta}(h)).$$ In other words, for any $\gamma\in H$, we have that $\gamma\in\mathfrak{A}_H(h)$ if and only if $\gamma\mathfrak{A}_{\theta}(h)\subseteq\mathfrak{A}_{\theta}(h)$.
\end{coro}
\begin{proof}
By Lemma \ref{lem:atheta} \ref{lem:atheta1}, for an element $\gamma\in H$, we have that $\gamma\cdot\mathfrak{p}_L^h\subseteq\mathfrak{p}_L^h$ if and only if $\gamma\cdot(\mathfrak{A}_{\theta}(h)\cdot\theta)\subseteq\mathfrak{A}_{\theta}(h)\cdot\theta$. Equivalently, $(\gamma\cdot\mathfrak{A}_{\theta}(h))\cdot\theta\subseteq\mathfrak{A}_{\theta}(h)\cdot\theta$. Since the map $h\in H\mapsto h\cdot\theta\in L$ is injective by definition of Hopf-Galois structure, the statement follows.
\end{proof}

We now establish some characterizations for the $\mathfrak{A}_H(h)$-freeness of $\mathfrak{p}_L^h$, which are the natural analogues of the ones for the case $h=0$ (see for instance \cite[Proposition 2.1]{gildegpgeneral}).

\begin{pro}\label{pro:hnormalbasisfract} Let $h\in\mathbb{Z}$. The following statements are equivalent:
\begin{enumerate}[label=\arabic*)]
    \item\label{pro:hnormalbasisfract1} $\mathfrak{p}_L^h$ is $\mathfrak{A}_H(h)$-free.
    \item\label{pro:hnormalbasisfract2} There is an $H$-normal basis generator $\theta\in\mathfrak{p}_L^h$ for $L$ such that $\mathfrak{A}_H(h)=\mathfrak{A}_{\theta}(h)$.
    \item\label{pro:hnormalbasisfract3} There is an $H$-normal basis generator $\theta\in\mathfrak{p}_L^h$ for $L$ such that $\mathfrak{A}_{\theta}(h)$ is a ring (namely, a subring of the underlying ring structure of $H$).
    \item\label{pro:hnormalbasisfract4} There is an $H$-normal basis generator $\theta\in\mathfrak{p}_L^h$ for $L$ such that $\mathfrak{A}_{\theta}(h)$ is principal as $\mathfrak{A}_H(h)$-fractional ideal. 
    \item\label{pro:hnormalbasisfract5} For every $H$-normal basis generator $\theta\in\mathfrak{p}_L^h$ for $L$, $\mathfrak{A}_{\theta}(h)$ is principal as $\mathfrak{A}_H(h)$-fractional ideal.
\end{enumerate}
\end{pro}
\begin{proof}
From the definition of $\mathfrak{A}_{\theta}(h)$, it is immediate that $\mathfrak{p}_L^h=\mathfrak{A}_{\theta}(h)\cdot\theta$. Thus, the equivalence between \ref{pro:hnormalbasisfract1} and \ref{pro:hnormalbasisfract2} is trivial. It is also immediate that \ref{pro:hnormalbasisfract2} implies \ref{pro:hnormalbasisfract3}. Conversely, suppose that for some $H$-normal basis generator $\theta\in\mathfrak{p}_L^h$, $\mathfrak{A}_{\theta}(h)$ is a ring. Given $\lambda\in\mathfrak{A}_{\theta}$, $$\lambda\mathfrak{p}_L^h=\lambda\mathfrak{A}_{\theta}\cdot\theta\subseteq\mathfrak{A}_{\theta}\cdot\theta=\mathfrak{p}_L^h.$$ We get $\lambda\in\mathfrak{A}_H(h)$, so $\mathfrak{A}_{\theta}(h)\subseteq\mathfrak{A}_H(h)$, whence $\mathfrak{A}_{\theta}(h)=\mathfrak{A}_H(h)$. Thus, \ref{pro:hnormalbasisfract2} follows. 

Next, we prove that \ref{pro:hnormalbasisfract1} implies \ref{pro:hnormalbasisfract5}. Assume that $\mathfrak{p}_L^h$ is $\mathfrak{A}_H(h)$-free, say $\mathfrak{p}_L^h=\mathfrak{A}_H(h)\cdot\alpha$ with $\alpha\in\mathfrak{p}_L^h$. Let $\theta\in\mathfrak{p}_L^h$ be an $H$-normal basis generator for $L$. Since $\mathfrak{A}_H(h)\cdot\alpha=\mathfrak{p}_L^h=\mathfrak{A}_{\theta}(h)\cdot\theta$, there is some $\varphi\in\mathfrak{A}_{\theta}(h)$ such that $\alpha=\varphi\cdot\theta$. Let $\lambda\in\mathfrak{A}_{\theta}(h)$. Then $\lambda\cdot\theta\in\mathfrak{p}_L^h=\mathfrak{A}_H(h)\cdot\alpha$, so there is $\gamma\in\mathfrak{A}_H(h)$ such that $$\lambda\cdot\theta=\gamma\cdot\alpha=\gamma\cdot(\varphi\cdot\theta)=(\gamma\cdot\varphi)\cdot\theta.$$ Since the map $j$ is bijective, $\lambda=\gamma\cdot\varphi\in\mathfrak{A}_H(h)\cdot\varphi$. We obtain that $\mathfrak{A}_{\theta}(h)=\mathfrak{A}_H(h)\cdot\varphi$, so \ref{pro:hnormalbasisfract5} holds.

It is trivial that \ref{pro:hnormalbasisfract5} implies \ref{pro:hnormalbasisfract4}. Finally, we prove that \ref{pro:hnormalbasisfract4} implies \ref{pro:hnormalbasisfract1}. Suppose that there is some $H$-normal basis generator $\theta\in\mathfrak{p}_L^h$ such that $\mathfrak{A}_{\theta}(h)$ is principal as $\mathfrak{A}_H(h)$-fractional ideal, say $\mathfrak{A}_{\theta}(h)=\mathfrak{A}_H(h)\cdot\varphi$. Then $$\mathfrak{p}_L^h=\mathfrak{A}_{\theta}(h)\cdot\theta=(\mathfrak{A}_H(h)\cdot\varphi)\cdot\theta=\mathfrak{A}_H(h)\cdot(\varphi\cdot\theta),$$ which proves that $\mathfrak{p}_L^h$ is $\mathfrak{A}_H(h)$-free.
\end{proof}

\section{Degree $p$ extensions of $p$-adic fields}\label{sec:degpsetup}

Let $p$ be an odd prime and let $L/K$ be a degree $p$ extension of $p$-adic fields. Let $\widetilde{L}$ be the normal closure of $L/K$ and $G=\mathrm{Gal}(\widetilde{L}/K)$. We know from a theorem by Galois that $G\cong C_p\rtimes C_r$ with $r\mid p-1$. Actually, $r=[\widetilde{L}:L]$. Fix a presentation $$G=\langle\sigma,\tau\,\mid\,\sigma^p=\tau^r=\mathrm{Id}_G,\,\tau\sigma=\sigma^g\tau\rangle,$$ where $g$ mod $p$ has order $r$ in $(\mathbb{Z}/p\mathbb{Z})^{\times}$. Call $J\coloneqq\langle\sigma\rangle$ and $G'=\langle\tau\rangle$. We can assume without loss of generality that $G'=\mathrm{Gal}(\widetilde{L}/L)$; otherwise we would replace $\tau$ by $\sigma^i\tau$ for a suitable $0\leq i\leq p-1$. Since $J$ is the only order $p$ subgroup of $G$, $M\coloneqq\widetilde{L}^J$ is the only degree $r$ intermediate field of $\widetilde{L}/K$. The Galois correspondence also yields that $L\cap M=K$ and $M/K$ Galois, so $L$ and $M$ are $K$-linearly disjoint. Then, the extension $M/K$ is cyclic of degree $r$. Since $p\nmid r$, by Hensel's lemma, $K$ contains the primitive $r$-th roots of unity. Hence, $M/K$ is cyclic and Kummer, so there is some $y\in M$ with $y^r\in K$ such that $M=K(y)$.

Call $s\coloneqq\frac{p-1}{r}$ and let $g_0\in\mathbb{Z}$ be such that $(\mathbb{Z}/p\mathbb{Z})^{\times}=\langle g_0\rangle$. Then $g_0^s$ and $g$ have the same order in $(\mathbb{Z}/p\mathbb{Z})^*$. Multiplying by a unit if necessary, one can assume that $g_0^s\equiv g\,(\mathrm{mod}\,p)$. Let $\chi\colon(\mathbb{Z}/p\mathbb{Z})^{\times}\longrightarrow\mathbb{Z}_p^{\times}$ be the $p$-adic multiplicative Teichmuller character. Let $\zeta_r$ be a primitive $r$-th root of unity such that $\chi(g)=\zeta_r$. Likewise, let $\zeta_{p-1}$ be a primitive $p-1$-th root of unity such that $\zeta_{p-1}=\chi(g_0)$. Then $\zeta_{p-1}^s=\zeta_r$.

Since $G\cong C_p\rtimes C_r$, we have that $G$ is solvable. Then, $L/K$ is Hopf-Galois and admits a unique Hopf-Galois structure $H$. We shall write $\mathfrak{A}_{L/K}$ (resp. $\mathfrak{A}_{L/K}(h)$) for the associated order of $\mathcal{O}_L$ (resp. $\mathfrak{p}_L^h$) in $H$. The main problem we consider in this paper is to characterize the $\mathfrak{A}_{L/K}(h)$-freeness of $\mathfrak{p}_L^h$ for each $h\in\mathbb{Z}$.

We know that the study of this kind of questions is closely related with the ramification on $L/K$. The only case in which $L/K$ can be possibly unramified is if it is Galois. A complete answer for the problem is easily found for this case. Indeed, assume that $L/K$ is unramified and let $G$ be its Galois group. The result \cite[Theorem 3.4]{truman2011} corresponds to an affirmative answer for $h=0$; namely, if $L/K$ is unramified, then $\mathfrak{A}_{L/K}=\mathcal{O}_K[G]$ and $\mathcal{O}_L$ is $\mathfrak{A}_{L/K}(h)$-free. Now, by Lemma \ref{lem:equivfreenesses}, $\mathfrak{A}_{L/K}(h)=\mathfrak{A}_{L/K}$ and $\mathfrak{p}_L^h$ is $\mathfrak{A}_{L/K}(h)$-free.

From now on, we will always assume that $L/K$ is ramified. Since it is of prime degree, it is automatically totally ramified. Since $e(L/K)=p$, from Lemma \ref{lem:equivfreenesses}, in order to study the $\mathfrak{A}_{L/K}(h)$-module structure of $\mathfrak{p}_L^h$, it is enough to consider $0\leq h\leq p-1$. Namely, the problem we consider is:

\begin{prob} Let $L/K$ be a ramified degree $p$ extension of $p$-adic fields and let $h\in\mathbb{Z}$ with $0\leq h\leq p-1$. Find a necessary and sufficient condition for $\mathfrak{p}_L^h$ being $\mathfrak{A}_{L/K}(h)$-free.
\end{prob}

Prior to the study of this problem, we shall view the structure of the underlying Hopf algebra at the Hopf-Galois structure $H$ on $L/K$.

\subsection{The Hopf-Galois structure on $L/K$}

The structure of the underlying Hopf algebra at $H$ was described by Childs \cite{childs1989} when $\mathrm{char}(K)=0$. Later, Elder \cite{elder2018} proved that this description is valid when $L/K$ is typical (that is, $\mathcal{O}_L$ is not generated by a $p$-th root of a uniformizer), with no restriction on the characteristic. We shall see that we can replay this description regardless of whether $L/K$ is typical or not. At the level of the arguments used, what follows in this part is merely expository.

It is immediate that $L/K$ is almost classically Galois with complement $M$. Then, the only Hopf-Galois structure on $L/K$ can be described as $H=M[J]^{G'}$ together with the Galois action. A $K$-basis of $M[J]$ is given by $$\{y^i\sigma^j\,\mid\,0\leq i\leq r-1,\,0\leq j\leq p-1\}.$$ Then, a $K$-basis of $H$ is formed by the non-zero sums of the elements at each orbit for the action of $G'$ on the basis above.

For $i=j=0$, we have $$y^0\sigma^0=\mathrm{Id},\quad\tau^k\cdot\mathrm{Id}=\mathrm{Id}.$$ Thus we have a single orbit with the identity.

For $i\neq0$ and $j=0$, we get $$y^i\sigma^0=y^i\mathrm{Id};\quad\tau^k\cdot y^i\mathrm{Id}=\zeta_r^{ki}y^i,\quad0\leq k\leq r-1.$$ The sum of the elements in this orbit is $$\sum_{k=0}^{r-1}\tau^k\cdot y^i\mathrm{Id}=\sum_{k=0}^{r-1}\zeta_r^{ki}y^i=0.$$ Since the elements of this orbit add up to $0$, we do not obtain basis element here.

Now, let $0\leq i\leq r-1$ and $1\leq j\leq p-1$. Since for all $0\leq k\leq r-1$ and $l\in\mathbb{Z}$ we have $\tau^k\cdot\sigma^{lg}=\sigma^{lg^{k+1}}$, we can assume without loss of generality that $j\notin\langle g\rangle$. A system of representatives for $(\mathbb{Z}/p\mathbb{Z})^{\times}-\langle g\rangle$ is $\{g_0^t\}_{t=0}^{s-1}$. Given $0\leq t\leq s-1$ and $0\leq k\leq r-1$, $$\tau^k\cdot y^i\sigma^{g_0^t}=\zeta_r^{ik}y^i\sigma^{g_0^{t+sk}}.$$ We calculate the sum of the elements in the orbit: $$\sum_{k=0}^{r-1}\tau^k\cdot y^i\sigma^{g_0^t}=y^i\sum_{k=0}^{r-1}\zeta_r^{ik}\sigma^{g_0^{t+sk}}.$$ For each $0\leq i\leq r-1$ and each $0\leq t\leq s-1$, define $$\theta(i,t)\coloneqq\sum_{k=0}^{r-1}\zeta_r^{ik}\sigma^{g_0^{t+sk}}=\sum_{k=0}^{r-1}\zeta_{p-1}^{isk}\sigma^{g_0^{t+sk}}.$$ Then $\sum_{k=0}^{r-1}\tau^k\cdot y^i\sigma^{g_0^t}=y^i\theta(i,t)$. Together with $\mathrm{Id}$, these form a $K$-basis of $H$. Namely, $$\{\mathrm{Id}\}\cup\{y^i\theta(i,t)\,\mid\,0\leq i\leq r-1,\,0\leq t\leq s-1\}.$$

We shall rewrite this basis in a more convenient way for our purposes. Let us extend the definition of $\theta(i,t)$ to $0\leq i\leq p-1$ in the obvious way. For each $0\leq b\leq s-1$, $$\theta(i+br,t)=\sum_{k=0}^{r-1}\zeta_r^{(i+br)k}\sigma^{g_0^{t+sk}}=\sum_{k=0}^{r-1}\zeta_r^{ik}\sigma^{g_0^{t+sk}}=\theta(i,t).$$ Given $0\leq i\leq p-1$, define $$\Theta(i)\coloneqq\sum_{t=0}^{s-1}\zeta_{p-1}^{it}\theta(i,t).$$ Then, \begin{equation*}
    \begin{split}
        \Theta(i)&=\sum_{t=0}^{s-1}\zeta_{p-1}^{it}\sum_{k=0}^{r-1}\zeta_{p-1}^{isk}\sigma^{g_0^{t+sk}}
        \\&=\sum_{t=0}^{s-1}\sum_{k=0}^{r-1}\zeta_{p-1}^{i(t+sk)}\sigma^{g_0^{t+sk}}
        \\&=\sum_{k=0}^{p-2}\zeta_{p-1}^{ik}\sigma^{g_0^k}.
    \end{split}
\end{equation*} 

Now, fix $0\leq i\leq r-1$. Given $0\leq b\leq s-1$, we have $$\Theta(i+br)=\sum_{t=0}^{s-1}\zeta_{p-1}^{(i+br)t}\theta(i+br,t)=\sum_{t=0}^{s-1}\zeta_{p-1}^{(i+br)t}\theta(i,t).$$ Hence, $$y^i\Theta(i+br)=\sum_{t=0}^{s-1}\zeta_{p-1}^{(i+br)t}y^i\theta(i,t).$$ Then, $$K\cdot\{y^i\Theta(i+br)\,\mid\,0\leq b\leq s-1\}\subseteq K\cdot\{y^i\theta(i,t)\,\mid\,0\leq t\leq s-1\}.$$ Conversely, given $0\leq t\leq s-1$, \begin{equation*}
    \begin{split}
        \sum_{b=0}^{s-1}\zeta_{p-1}^{-t(i+br)}\Theta(i+br)&=\sum_{b=0}^{s-1}\zeta_{p-1}^{-t(i+br)}\sum_{u=0}^{s-1}\zeta_{p-1}^{u(i+br)}\theta(i,u)\\&=\sum_{b,u=0}^{s-1}\zeta_{p-1}^{(u-t)(i+br)}\theta(i,u)\\&=s\theta(i,t),
    \end{split}
\end{equation*} whence $y^i\theta(i,t)=\frac{1}{s}\sum_{b=0}^{s-1}\zeta_{p-1}^{-t(i+br)}y^i\Theta(i+br)$. We have proved that $$K\cdot\{y^i\Theta(i+br)\,\mid\,0\leq b\leq s-1\}=K\cdot\{y^i\theta(i,t)\,\mid\,0\leq t\leq s-1\}.$$ 

We get that $$\{\mathrm{Id}\}\cup\{y^i\Theta(i+br)\,\mid\,0\leq i\leq r-1,\,0\leq b\leq s-1\}$$ is a $K$-basis of $H$. Now, $y^{br}=(y^r)^b\in K^{\times}$, so $$\{\mathrm{Id}\}\cup\{y^{i+br}\Theta(i+br)\,\mid\,0\leq i\leq r-1,\,0\leq b\leq s-1\}$$ is also a $K$-basis of $H$. Since $i+br$ runs through $\{0,\dots,p-2\}$ as $i$ (resp. $b$) runs through $\{0,\dots,r-1\}$ (resp. $\{0,\dots,s-1\}$), we get that $$\{\mathrm{Id}\}\cup\{y^i\Theta(i)\,\mid\,0\leq i\leq p-2\}$$ is a $K$-basis of $H$. Recall that $$y^i\Theta(i)=y^i\sum_{k=0}^{p-2}\zeta_{p-1}^{ik}\sigma^{g_0^k}.$$

Given $1\leq i\leq p-2$, define $$\Psi_i\coloneqq-\frac{1}{y^i}\Theta(p-1-i).$$ Since $y^{p-1}=y^{rs}$, we can rewrite $$\Psi_i=-\frac{1}{(y^r)^s}y^{p-1-i}\Theta(p-1-i).$$ Since $y^r\in K$ and $y^{p-1-i}\Theta(p-1-i)\in H$, we have that $\Psi_i\in H$. What is more, we deduce that $\{\Psi_i\}_{i=1}^{p-2}$ is a $K$-linearly independent system. We prefer an alternative expression for these. Namely, we have $$\Psi_i=-\frac{1}{y^i}\sum_{k=0}^{p-2}\zeta_{p-1}^{-ik}\sigma^{g_0^k}.$$ Since $g_0$ is a generator of $(\mathbb{Z}/p\mathbb{Z})^{\times}$ and $\zeta_{p-1}=\chi(g_0)$, we get $$\Psi_i=-\frac{1}{y^i}\sum_{m=1}^{p-1}\chi(m)^{-i}\sigma^m.$$

Define also $\Psi_0=\mathrm{Id}$ and $$\Psi_{p-1}=\frac{p-1-\Theta(0)}{y^{p-1}}=\frac{p-1-\Big(\sum_{m=1}^{p-1}\sigma^m\Big)}{y^{p-1}}.$$ Then, we obtain that $\{\Psi_i\}_{i=0}^{p-1}$ is a $K$-basis of $H$.

Let us call $\Psi\coloneqq\Psi_1$. Then, it is proved that this element is a generator of $H$ as a $K$-algebra. Namely:

\begin{thm}\label{thm:hopfalggen} Let $L/K$ be a degree $p$ extension of $p$-adic fields with normal closure $\widetilde{L}$ and let $r=[\widetilde{L}:L]$. Let $M$ be the complement of $L$ inside $\widetilde{L}$, which is a cyclic Kummer degree $r$ extension of $K$. Let $y\in M$ be a Kummer generator with $\tau(y)=\zeta_ry$. Then $H=K[\Psi]$, where \begin{equation}\label{eq:hopfalggen}
    \Psi=-\frac{1}{y}\Big(\sum_{m=1}^{p-1}\chi(m)^{-1}\sigma^m\Big).
\end{equation} Moreover, there is a unit $\varepsilon\in\mathbb{Z}_p^{\times}$ such that the element $$\lambda\coloneqq\frac{\varepsilon p}{y^{p-1}}$$ satisfies $\Psi^p=\lambda\Psi$.
\end{thm}

This result is proved from the classification of rank $p$ Hopf algebras by Tate and Oort \cite{tateoort} consisting in eigenspace decompositions of $H$, as presented in \cite[Chapter 4]{childs}. We proceed to sketch the procedure so as to prove it (see \cite[Section 5.1.1]{gildegpgeneral} for a detailed explanation). Let $[\cdot]\colon\mathbb{Z}/p\mathbb{Z}\longrightarrow\mathrm{End}_{K-Hopf}(H)$ be the ring homomorphism from \cite[\textsection 16]{childs}. Let $H^+=\mathrm{Ker}(\varepsilon_H)$. For each $1\leq i\leq p-1$, call $$H_i=\{h\in H^+\,\mid\,[m](h)=\chi^i(m)h\hbox{ for all }1\leq m\leq p-1\}.$$ Then $H=\bigoplus_{i=0}^{p-1}H_i$ has the structure of a graded ring, that is, $H_iH_j\subset H_{i+j\mod p}$ for all $0\leq i,j\leq p-1$. Now, it is proved that $H_i=K\Psi_i$ for all $1\leq i\leq p-1$. Thus, if we call $\Psi_1\coloneqq\Psi$, then $\Psi^r\in H_r$ for all $1\leq r\leq p-1$, whence $H=K[\Psi]$. Moreover, $\Psi^p\in H_1$, whence $\Psi^p=\lambda\Psi$ with $\lambda\in K$. The description of $\lambda$ follows from \cite[(16.8)]{childs}.

Let $e\coloneqq e(K/\mathbb{Q}_p)$. Recall that $y\in M$ is any element such that $M=K(y)$ and $y^r\in K$. We have that $$v_K(\lambda)=e-\frac{p-1}{r}v_K(y^r).$$ In the case that $L/K$ is not maximally ramified, the extension is typical and admits a scaffold (this is proved in \cite[\textsection 3.2]{elder2018}). After a suitable choice of $y$ in such a way that $\lambda\in\mathcal{O}_K$, $v_K(\lambda)$ affects directly to the precision of the scaffold (see Proposition \ref{pro:precscaffold} below). Thus, this value is of great importance for the study of $\mathfrak{p}_L^h$ as an $\mathfrak{A}_H(h)$-module.

\section{Reduction to the totally ramified case}\label{sec:redtotram}

Let $L/K$ be a ramified degree $p$ extension of $p$-adic fields with normal closure $\widetilde{L}$, and acquire the notations from Section \ref{sec:degpsetup}. We have that $\widetilde{L}/K$ is not in general totally ramified. In fact, $$e(\widetilde{L}/K)=pe(\widetilde{L}/L)$$ with $e(\widetilde{L}/L)\mid r$. Since $L$ and $M$ are $K$-linearly disjoint, $e(\widetilde{L}/L)=e(M/K)$. Thus, $\widetilde{L}/K$ is totally ramified if and only if so is $\widetilde{L}/L$; if and only if so is $M/K$. In this part we will prove that, in order to characterize the $\mathfrak{A}_{L/K}(h)$-freeness of $\mathfrak{p}_L^h$, we can assume without loss of generality that $\widetilde{L}/K$ is totally ramified. In \cite[Section 3]{gildegpgeneral}, this was shown for the case $h=0$. We will see that we can recycle most of the work therein. The concrete result for that case is the following:

\begin{pro}[\cite{gildegpgeneral}, Proposition 3.4]\label{pro:redtotramh=0} Let $L/K$ be a ramified extension of degree $p$ with normal closure $\widetilde{L}$ and let $M$ be the complement of $L/K$. Let $L'$ (resp. $K'$) be the inertia field of $\widetilde{L}/L$ (resp. $\widetilde{L}/K$). Then:
\begin{enumerate}[label=\arabic*)]
    \item $L'/K'$ is a degree $p$ extension of $p$-adic fields with normal closure $\widetilde{L}$.
    \item $\widetilde{L}/K'$ is totally ramified.
    \item $L'/K'$ is almost classically Galois with complement $M$.
    \item $H'\coloneqq H\otimes_KK'$ is the only Hopf-Galois structure on $L'/K'$.
    \item\label{pro:redtotramh=05} $\mathfrak{A}_{L'/K'}=\mathcal{O}_{K'}\otimes_{\mathcal{O}_K}\mathfrak{A}_{L/K}$ and $\mathcal{O}_{L'}$ is $\mathfrak{A}_{L'/K'}$-free if and only if $\mathcal{O}_L$ is $\mathfrak{A}_{L/K}$-free.
\end{enumerate}
\end{pro}

It holds that $L'=LK'$ with $L$ and $K'$ $K$-linearly disjoint. The lattice with the involved intermediate fields of $\widetilde{L}/K$ is as follows: 

\[
\xymatrix{
&& \ar@{-}@/_2pc/[ddll]_{G'} \ar@{-}[dl]_{G_0'} \widetilde{L} \ar@{-}[dd]^{G_0} \ar@{-}[dr]^J & \\
& \ar@{-}[dl] L' \ar@{-}[dr]_{H'} & & \ar@{-}[dl] M \\
L \ar@{-}[dr]_H & & \ar@{-}[dl] K' & \\
& K & &
}
\]

We now prove the natural generalization of Proposition \ref{pro:redtotramh=0} \ref{pro:redtotramh=05} to the ideals $\mathfrak{p}_L^h$, $0\leq h\leq p-1$.

\begin{pro}\label{pro:redtotramh} Take all the notations from Proposition \ref{pro:redtotramh=0}. Given $0\leq h\leq p-1$, we have that $\mathfrak{A}_{L'/K'}(h)=\mathfrak{A}_{L/K}(h)\otimes_{\mathcal{O}_K}\mathcal{O}_{K'}$, and $\mathfrak{p}_{L'}^{h}$ is $\mathfrak{A}_{L'/K'}(h)$-free if and only if $\mathfrak{p}_L^h$ is $\mathfrak{A}_{L/K}(h)$-free.
\end{pro}
\begin{proof}
First, let us prove that $\mathfrak{A}_{L'/K'}(h)=\mathfrak{A}_{L/K}(h)\otimes_{\mathcal{O}_K}\mathcal{O}_{K'}$. Since $K'/K$ is unramified, it is arithmetically disjoint with $L/K$. Then, $\mathcal{O}_{L'}=\mathcal{O}_L\otimes_{\mathcal{O}_K}\mathcal{O}_{K'}$. Also, since $L'/L$ is unramified, a uniformizer of $L$ is also a uniformizer of $L'$, so $\mathfrak{p}_{L'}^h=\mathfrak{p}_L^h\otimes_{\mathcal{O}_K}\mathcal{O}_{K'}$. Thus, $\mathfrak{A}_{L/K}(h)\otimes_{\mathcal{O}_K}\mathcal{O}_{K'}$ acts componentwise on $\mathfrak{p}_{L'}^h$ through the factors of the tensor products, and hence $\mathfrak{A}_{L/K}(h)\otimes_{\mathcal{O}_K}\mathcal{O}_{K'}\subseteq\mathfrak{A}_{L'/K'}(h)$. Conversely, let $\gamma\in\mathfrak{A}_{L'/K'}(h)$, so that $\gamma\in H'=H\otimes_KK'$. Let $\{z_j\}_{j=1}^s$ be an $\mathcal{O}_K$-basis of $\mathcal{O}_{K'}$. Then, there are unique elements $\gamma_1,\dots,\gamma_s\in H$ such that $$\gamma=\sum_{j=1}^s\gamma_jz_j.$$ It is enough to prove that $\gamma_j\in\mathfrak{A}_{L/K}(h)$ for all $1\leq j\leq s$. Indeed, given $\theta\in\mathfrak{p}_L^h\subseteq\mathfrak{p}_{L'}^h$, we have that $\gamma\cdot\theta\in\mathfrak{p}_{L'}^h$. Now, $$\gamma\cdot\theta=\sum_{j=1}^s(\gamma_j\cdot\theta)z_j.$$ Since $\mathfrak{p}_L^h$ is $\mathcal{O}_K$-flat, $\{z_j\}_{j=1}^s$ is a $\mathfrak{p}_L^h$-basis of $\mathfrak{p}_{L'}^h$. Therefore, it follows that $\gamma_j\cdot\theta\in\mathfrak{p}_L^h$ for all $1\leq j\leq s$, proving that $\gamma_j\in\mathfrak{A}_{L/K}(h)$ as desired.

Let us prove that $\mathfrak{p}_{L'}^{h}$ is $\mathfrak{A}_{L'/K'}(h)$-free if and only if $\mathfrak{p}_L^h$ is $\mathfrak{A}_{L/K}(h)$-free. The right to left implication follows from the fact that $\mathcal{O}_{K'}$ is $\mathcal{O}_K$-flat. Conversely, assume that $\mathfrak{p}_{L'}^{h}$ is $\mathfrak{A}_{L'/K'}(h)$-free. Since $K'/K$ is unramified, $\mathcal{O}_{K'}$ is $\mathcal{O}_K$-faithfully flat. Now, we use faitfully flat descent. From the hypothesis, $\mathfrak{p}_{L'}^{h}$ is $\mathfrak{A}_{L'/K'}(h)$-projective of rank one, so $\mathfrak{p}_L^{h}$ is $\mathfrak{A}_{L/K}(h)$-projective of rank one. Now, $\mathfrak{A}_{L/K}(h)$ is a commutative semilocal ring, so every projective module of rank one over $\mathfrak{A}_{L/K}(h)$ is free. We conclude that $\mathfrak{p}_L^h$ is $\mathfrak{A}_{L/K}(h)$-free, as desired.
\end{proof}

\section{The maximally ramified case}\label{sec:maxram}

Let $L/K$ be a ramified degree $p$ extension of $p$-adic fields with normal closure $\widetilde{L}$ and let $0\leq h\leq p-1$. By Proposition \ref{pro:redtotramh}, there is no loss of generality in assuming that $\widetilde{L}/K$ is totally ramified so as to characterize the $\mathfrak{A}_{L/K}(h)$-freeness of $\mathrm{p}_L^h$. Let $t$ be the ramification jump of $\widetilde{L}/K$. Then, it is checked that the ramification jump of $L/K$ is $\ell=\frac{t}{r}$, where $r=[\widetilde{L}:L]$ (see \cite[Section 2.2.6]{elder2018}). Note that $\ell$ is a $p$-adic integer that is not in general a rational integer. The number $t$ is also the ramification jump of $\widetilde{L}/M$, which is a cyclic degree $p$ extension of $p$-adic fields. The following result is contained in \cite[Chapter III, \textsection2, Proposition 2.3]{fesenkovostokov}.

\begin{pro}\label{pro:ramifdegp} Let $E/F$ be a ramified degree $p$ Galois extension of $p$-adic fields with ramification jump $t$. Then, $$1\leq t\leq\frac{pe(F/\mathbb{Q}_p)}{p-1}.$$ If in addition $p\mid t$, then $t=\frac{pe(F/\mathbb{Q}_p)}{p-1}$, $F$ contains the $p$-th roots of unity and $E=F(\gamma)$ for some $\gamma\in E$ with $\gamma^p\in F$ and $v_E(\gamma)=1$.
\end{pro}

We have actually that $t=\frac{pe(F/\mathbb{Q}_p)}{p-1}$ if and only if $p\mid t$ (see \cite[(24.7)]{childs}).

In our case, the extension $\widetilde{L}/M$ is under the hypotheses of Proposition \ref{pro:ramifdegp}. Note that its ramification jump is $t$ (because of \cite[Chapter IV, \textsection1, Proposition 2]{serre}) and, calling $e=e(K/\mathbb{Q}_p)$, $e(M/\mathbb{Q}_p)=re$. We obtain that $1\leq t\leq\frac{rpe}{p-1}$. Now, since $\ell=\frac{t}{r}$, we have \begin{equation}
    0<\ell\leq\frac{pe}{p-1}.
\end{equation}

In this section, we will consider the case in which $\widetilde{L}/M$ is maximally ramified, which means that $t=\frac{rpe}{p-1}$. Equivalently, $\ell=\frac{pe}{p-1}$. Consequently, we will also say that $L/K$ is maximally ramified. As already mentioned, all of this is equivalent to $p\mid t$. Our main aim is to prove that $\mathfrak{p}_L^h$ is $\mathfrak{A}_{L/K}(h)$-free for all $0\leq h\leq p-1$.

The case for $h=0$ was considered in \cite[Section 4.1]{gildegpgeneral}, where it was proved that $\mathfrak{A}_{L/K}$ is the maximal $\mathcal{O}_K$-order in $H$. We shall prove that $\mathfrak{A}_{L/K}(h)=\mathfrak{A}_{L/K}$ for all $0\leq h\leq p-1$, whence the maximality of $\mathfrak{A}_{L/K}(h)$ will follow. By \cite[Proposition 2.5.5]{trumanthesis}, this automatically implies that $\mathfrak{p}_L^h$ is $\mathfrak{A}_{L/K}(h)$-free for all $0\leq h\leq p-1$.

From Proposition \ref{pro:ramifdegp}, we know that $\widetilde{L}=M(\gamma)$ with $\gamma^p\in M$ and $v_{\widetilde{L}}(\gamma)=1$, and $M$ contains the $p$-th roots of unity. In other words, $\widetilde{L}/M$ is a cyclic degree $p$ Kummer extension with Kummer generator $\gamma$. Recall that the Galois group of $\widetilde{L}/M$ is $J=\langle\sigma\rangle$. Then, for a suitable choice of primitive $p$-th root of unity $\xi_p$, we have that $\sigma(\gamma)=\xi_p\gamma$. We will prefer to work with $\zeta_p\coloneqq\xi_p^r$, which is also a primitive $p$-th root of unity because $r\mid p-1$ implies $p\nmid r$.

\begin{pro}\label{pro:galmaxram} Let $\alpha=N_{\widetilde{L}/L}(\gamma)\in L$. Then:
\begin{enumerate}[label=\arabic*)]
    \item\label{pro:galmaxram1} $L=K(\alpha)$ with $\alpha^p\in K$ and $v_L(\alpha)=1$.
    \item\label{pro:galmaxram2} $M=K(\zeta_p)$.
    \item\label{pro:galmaxram3} The action of $G$ on $\widetilde{L}$ is given by $$\sigma(\alpha)=\zeta_p\alpha,\quad\sigma(\zeta_p)=\zeta_p,$$ $$\tau(\alpha)=\alpha,\quad\tau(\zeta_p)=\zeta_p^g.$$
\end{enumerate}
\end{pro}
\begin{proof}
The statements \ref{pro:galmaxram1} and \ref{pro:galmaxram2} are proved in \cite[Proof of Lemma 7.2]{gilkummer}. We deduce that $\widetilde{L}=K(\alpha,\zeta_p)$, so the action of $G$ on $\widetilde{L}$ is completely determined by the action of $\sigma$ and $\tau$ on $\alpha$ and $\zeta_p$. Since $\alpha\in L$ and $\zeta_p\in M$, we have $\sigma(\zeta_p)=\zeta_p$ and $\tau(\alpha)=\alpha$ (under the convention that $\tau$ is the order $r$ element of $G$ that fixes $L$). Next, let us prove that $\sigma(\alpha)=\zeta_p\alpha$. Since $G'=\langle\tau\rangle$, $$\alpha=N_{\widetilde{L}/L}(\gamma)=\prod_{k=0}^{r-1}\tau^k(\gamma).$$ Now, \begin{equation*}
    \begin{split}
        \sigma\tau^k(\gamma)&=\tau^k\sigma^{g^{-k}}(\gamma)
        \\&=\tau^k(\xi_p^{g^{-k}}\gamma)
        \\&=\xi_p\tau^k(\gamma)
    \end{split}
\end{equation*} Hence, \begin{equation*}
    \begin{split}
        \sigma(\alpha)&=\prod_{k=0}^{r-1}\sigma\tau^k(\gamma)
        \\&=\prod_{k=0}^{r-1}(\xi_p\tau^k(\gamma))
        \\&=\xi_p^r\prod_{k=0}^{r-1}\tau^k(\gamma)
        \\&=\zeta_p\alpha,
    \end{split}
\end{equation*}
proving the formula for $\sigma$. Finally, $$\tau(\zeta_p)=\frac{\tau\sigma(\alpha)}{\alpha}=\frac{\zeta_p^g\alpha}{\alpha}=\zeta_p^g,$$ so \ref{pro:galmaxram3} follows.
\end{proof}

This simple description for the Galois action of $\widetilde{L}/M$ can be used to prove the maximality of $\mathfrak{A}_{L/K}$ through the idempotents of $H$. First, we state a description for the idempotents of a Galois group algebra.

\begin{pro}\label{pro:idempotgroupalg} Let $F$ be a field with $\mathrm{char}(F)=0$ and let $\Gamma$ be a finite abelian group. Let $I$ be the set of irreducible $F$-characters of $\Gamma$. The primitive idempotents of $F[\Gamma]$ are of the form $$e_{\chi}=\frac{1}{|\Gamma|}\sum_{\gamma\in\Gamma}\chi(\gamma^{-1})\gamma,\quad\chi\in I.$$ What is more, $\{e_{\chi}\}_{\chi\in I}$ is a primitive pairwise-orthogonal system of idempotents for $F[\Gamma]$.
\end{pro}

These are standard facts in representation theory of finite groups and can be found at \cite[\textsection1.4]{bertrandias1977}. In our case, we apply Proposition \ref{pro:idempotgroupalg} to $M[J]$. Since $M=K(\zeta_p)$, the irreducible $M$-characters of $J$ are $$\chi_i(\sigma)=\zeta_p^i,\quad0\leq i\leq p-1.$$ For each such an $i$, let us call $e_i\coloneqq e_{\chi_i}$.

\begin{coro} The primitive idempotents of $M[J]$ are \begin{equation}\label{eq:idempot}
    e_i=\frac{1}{p}\sum_{m=0}^{p-1}\zeta_p^{-im}\sigma^m\in M[J],\quad0\leq i\leq p-1.
\end{equation} In addition, $\{e_i\}_{i=0}^{p-1}$ is a primitive pairwise orthogonal system of idempotents for $M[J]$.
\end{coro}

From here, one can prove the maximality of $\mathfrak{A}_{\widetilde{L}/M}$ as an $\mathcal{O}_M$-order: because $e_0,\dots,e_{p-1}$ form a system of primitive pairwise orthogonal idempotents adding up to $1$, the map $M[J]\longrightarrow M^p$ defined by $e_i\mapsto(\delta_{ij})_{j=0}^{p-1}$ is an isomorphism of $M$-algebras that sends $\mathfrak{A}_{\widetilde{L}/M}$ to $\mathcal{O}_M^p$. Since the latter is the maximal $\mathcal{O}_M$-order of $M^p$, the former is the maximal $\mathcal{O}_M$-order of $\mathfrak{A}_{\widetilde{L}/M}$.

In order to prove the analogous result for $L/K$, we descend the data from $\widetilde{L}/M$ to $L/K$. We already saw that $L=K(\alpha)$ with $\alpha^p\in K$ and $v_L(\alpha)=1$. This last equality says that $\alpha$ is a uniformizer of $L$, so $\{\alpha^j\}_{j=0}^{p-1}$ is an integral basis of $L/K$. On the other hand, the idempotents $e_0,\dots,e_{p-1}$ descend without any change.

\begin{pro} For $0\leq i\leq p-1$, let $e_i$ be as in \eqref{eq:idempot}. Then, $\{e_i\}_{i=0}^{p-1}$ is a primitive pairwise orthogonal system of idempotents for $H$.
\end{pro}
\begin{proof}
We have that \begin{equation*}
    \begin{split}
        \tau*e_i&=\frac{1}{p}\sum_{m=0}^{p-1}\zeta_p^{-igm}\sigma^{gm}
        \\&=\frac{1}{p}\sum_{m=0}^{p-1}\zeta_p^{-im}\sigma^m
        \\&=e_i.
    \end{split}
\end{equation*} Hence, $e_i\in M[J]^{G'}=H$ for all $0\leq i\leq p-1$. Since in addition $M[J]=H\otimes_KM$ and $\{e_i\}_{i=0}^{p-1}$ is a primitive pairwise orthogonal system for $M[J]$, so is for $H$.
\end{proof}

Now, we can repeat the argument as before so as to prove the maximality of $\mathfrak{A}_{L/K}$ in $H$.

\begin{pro}[\cite{gildegpgeneral}, Corollary 4.2]\label{pro:maximalassocordh=0} The associated order $\mathfrak{A}_{L/K}$ of $\mathcal{O}_L$ in $H$ is the maximal $\mathcal{O}_K$-order in $H$.
\end{pro}
\begin{proof}
For $0\leq i\leq p-1$, let $e_i$ be as in \eqref{eq:idempot}. Since $\sigma(\gamma)=\zeta_p^{r^{-1}}\gamma$ and $\alpha=N_{\widetilde{L}/L}(\gamma)$, we have that $\sigma(\alpha)=\zeta_p\alpha$. Now, given $0\leq i,j\leq p-1$, \begin{equation}\label{eq:diagaction}
    \begin{split}
        e_i\cdot\alpha^j&=\frac{1}{p}\sum_{m=0}^{p-1}\zeta_p^{-im}\sigma^m(\alpha^j)
        \\&=\frac{1}{p}\sum_{m=0}^{p-1}\zeta_p^{-im}\zeta_p^{mj}\alpha^j
        \\&=\frac{1}{p}\sum_{m=0}^{p-1}\zeta_p^{(j-i)m}\alpha^j
        \\&=\delta_{ij}\alpha^j
    \end{split}
\end{equation} We see that the idempotents $e_i$ act diagonally on the powers of $\alpha$. Since $\{\alpha^j\}_{j=0}^{p-1}$ is an $\mathcal{O}_K$-basis of $\mathcal{O}_L$, it follows that $\{e_i\}_{i=0}^{p-1}$ is an $\mathcal{O}_K$-basis of $\mathfrak{A}_{L/K}$. Then, the map $H\longrightarrow K^p$ defined by $e_i\mapsto(\delta_{ij})_{j=0}^{p-1}$ is an isomorphism of $K$-algebras sending $\mathfrak{A}_{L/K}$ to $\mathcal{O}_K^p$. Since the latter is the maximal $\mathcal{O}_K$-order in $K^p$, $\mathfrak{A}_{L/K}$ is the maximal $\mathcal{O}_K$-order in $H$.
\end{proof}

\begin{rmk}\normalfont The contents in this section up to Proposition \ref{pro:maximalassocordh=0} are essentially expository, and have culminated in a proof of the maximality of $\mathfrak{A}_{L/K}$, by finding a system of primitive pairwise orthogonal idempotents of $H$, namely $\{e_i\}_{i=0}^{p-1}$. In \cite[\textsection4]{gilkummer}, our elements $e_i$ are denoted $v_i$, and they are found from the action of the $K$-basis $\{\Psi^i\}_{i=0}^{p-1}$ on the integral basis $\{\alpha^j\}_{j=0}^{p-1}$ of $L/K$ by using \cite[Theorem 1.2]{gilkummer} (see \cite[\textsection 4]{gildegpgeneral} and the preceding discussion). At the same time, that last result makes use of a more general method from Rio and the author \cite{gilrioinduced}, which studies the Hopf-Galois module structure of $L/K$ consisting in the study of the action matrix of a basis of $H$ on an integral basis for $L/K$. Actually, there is no necessity for the use of this method because the elements $e_i$ are already known from the Galois case.
\end{rmk}


Now, we consider the module structure of the ideals $\mathfrak{p}_L^h$, $0\leq h\leq p-1$. Since $\{\alpha^j\}_{j=0}^{p-1}$ is an $\mathcal{O}_K$-basis of $\mathcal{O}_L$, we deduce immediately that $\{\alpha^{h+j}\}_{j=0}^{p-1}$ is an $\mathcal{O}_K$-basis of $\mathfrak{p}_L^h$. 

\begin{pro}\label{pro:maximalassocordgenh} Let $L/K$ be a maximally ramified degree $p$ extension of $p$-adic fields and let $H$ be its only Hopf-Galois structure. Given $0\leq h\leq p-1$, $\mathfrak{A}_{L/K}(h)=\mathfrak{A}_{L/K}$ and $\mathfrak{p}_L^h$ is $\mathfrak{A}_{L/K}(h)$-free. Moreover, the free generators of $\mathfrak{p}_L^h$ as an $\mathfrak{A}_{L/K}(h)$-module are $$\beta=\sum_{j=0}^{p-1}\beta_{h+j}\alpha^{h+j},\quad\beta_{h+j}\in\mathcal{O}_K^{\times}.$$
\end{pro}
\begin{proof}
Fix $0\leq h\leq p-1$. For each $0\leq i,j\leq p-1$, we determine $e_i\cdot\alpha^{h+j}$ using \eqref{eq:diagaction}. If $j\leq p-h-1$ (equivalently, $h+j\leq p-1$), we have $$e_i\cdot\alpha^{h+j}=\delta_{i,h+j}\alpha^{h+j}.$$ On the contrary, if $j\geq p-h$, then $h+j-p<p$. Since $\alpha^p\in K$, \begin{equation*}
    \begin{split}
        e_i\cdot\alpha^{h+j}&=\alpha^pe_i\cdot\alpha^{h+j-p}
        \\&=\alpha^p\delta_{i,h+j-p}\alpha^{h+j-p}
        \\&=\delta_{i,h+j-p}\alpha^{h+j}
    \end{split}
\end{equation*} All together, $$e_i\cdot\alpha^{h+j}=\delta_{i,[h+j]_p}\alpha^{h+j},\quad 0\leq i,j\leq p-1.$$ Note that $j\mapsto[h+j]_p$ defines a permutation of $\{0,\dots,p-1\}$. Since $\{\alpha^{h+j}\}_{j=0}^{p-1}$ is an $\mathcal{O}_K$-basis of $\mathfrak{p}_L^h$, we find that $\{e_i\}_{i=0}^{p-1}$ is an $\mathcal{O}_K$-basis for $\mathfrak{A}_{L/K}(h)$. Since it also was for $\mathfrak{A}_{L/K}$, the equality $\mathfrak{A}_{L/K}=\mathfrak{A}_{L/K}(h)$ follows. As already mentioned, this equality implies the $\mathfrak{A}_{L/K}(h)$-freeness of $\mathfrak{p}_L^h$.

Let $\beta\in\mathfrak{p}_L^h$, and write $\beta=\sum_{j=0}^{p-1}\beta_{h+j}\alpha^{h+j}$ with $\beta_{h+j}\in\mathcal{O}_K$. Then, for each $0\leq i\leq p-1$.  \begin{equation*}
    \begin{split}
        e_i\cdot\beta&=\sum_{j=0}^{p-1}\beta_{h+j}e_i\cdot\alpha^{h+j}
        \\&=\sum_{j=0}^{p-h-1}\beta_{h+j}\delta_{i,h+j}\alpha^{h+j}+\sum_{j=p-h}^{p-1}\beta_{h+j}\delta_{i,h+j-p}\alpha^{h+j},
    \end{split}
\end{equation*} the second sum appearing unless $h=0$. If $h\leq i\leq p-1$, then the only non-vanishing summand of $e_i\cdot\beta$ is at the first sum, namely for $i=h+j$ with $0\leq j\leq p-1$. Hence, $e_i\cdot\beta=\beta_i\alpha^i$. Otherwise, if $0\leq i<h$, then the only non-vanishing summand is at the second sum, for $i=h+j-p$ with $p-h\leq j\leq p-1$. Thus, $e_i\cdot\beta=\beta_{p+i}\alpha^{p+i}$. Then, $$\{e_i\cdot\beta\}_{i=0}^{p-1}=\{\beta_{h+j}\alpha^{h+j}\}_{j=0}^{p-1}.$$ We conclude that $\beta$ is an $\mathfrak{A}_{L/K}(h)$-free generator if and only if $\beta_{h+j}\in\mathcal{O}_K^{\times}$ for all $0\leq j\leq p-1$.
\end{proof}

\section{The case of typical degree $p$ extensions}\label{sec:typicaldegp}

In this section, we consider the case in which our (ramified) degree $p$ extension $L/K$ of $p$-adic fields is not maximally ramified. Namely, the ramification jump $\ell$ of $L/K$ satisfies $$0<\ell<\frac{pe}{p-1}.$$ Recall that this holds if and only if $p$ does not divide the ramification jump $t$ of $\widetilde{L}/K$. Equivalently, $\ell\not\equiv0\,(\mathrm{mod}\,p)$ in the ring $\mathbb{Z}_p$ of $p$-adic integers. Let $a$ be the residue class of $\ell$ mod $p$. Then, this is the case where $a\neq0$.

This kind of extensions was extensively studied by Elder in \cite{elder2018}, where he refers to them as \textbf{typical}. Namely, a separable ramified extension of local fields $E/F$ with degree $p$ and residue characteristic $p$ is said to be typical if it is not generated by a $p$-th root of a uniformizer $\pi_F$. In our case, $L/K$ is typical if and only if it is not maximally ramified (see for instance \cite[(24.2) and (24.3)]{childs}).

In \cite[Theorem 2.2]{elder2018}, Elder provides a characterization of a typical extension in terms of the equation satisfied by some primitive element of the extension. We state here a summary of that result for the case in which $\widetilde{L}/K$ is totally ramified.

\begin{thm}\label{thm:eldertypicalp} Let $L/K$ be a degree $p$ extension of $p$-adic fields with normal closure $\widetilde{L}$. Assume that $\widetilde{L}/K$ is totally ramified and call $r=[\widetilde{L}:L]$. The extension $L/K$ is typical if and only if there is a primitive element $x$ of $L/K$, $\alpha,\beta\in K$ and integers $b,c\in\mathbb{Z}$ with $0\leq c<r$ and $\gcd(c,r)=1$ such that \begin{equation}\label{eq:typicaldegp}
    x^p-\alpha^{\frac{p-1}{r}}x-\beta=0,
\end{equation} $v_K(\alpha)=c$ and $v_K(\beta)=-b$. In that case, the ramification jump of $L/K$ is \begin{equation}\label{eq:fractramjump}
    \ell=b+\frac{pc}{r}.
\end{equation}
\end{thm}

This statement provides a characterization of whether a given degree $p$ extension $L/K$ of $p$-adic fields is typical. But actually, Elder's result is stronger: if $b,c,t,r\in\mathbb{Z}$ are in the conditions of Theorem \ref{thm:eldertypicalp}, then there is a typical degree $p$ extension $L$ of $K$ with a primitive element $x$ satisfying \eqref{eq:typicaldegp} and with ramification jump as in \eqref{eq:fractramjump} (see also the statement \cite[Theorem 5.1]{gildegpgeneral}).

Let $b,c\in\mathbb{Z}$ and $\alpha,\beta\in K$ be as in Theorem \ref{thm:eldertypicalp} for our extension $L/K$. From the equation \eqref{eq:fractramjump} we see that $\ell\equiv b\,(\mathrm{mod}\,p)$ in $\mathbb{Z}_p$, whence $a$ is the remainder of $b$ mod $p$. On the other hand, recall that, for the complement $M$ of $L/K$, $M=K(y)$ with $y^r\in K$. In particular, one can choose $y$ such that $y^r=\alpha$, so that $v_M(y)=c$. This will be our choice for $y$ from now on, also in the expression of the element $\Psi$ such that $H=K[\Psi]$ from \eqref{eq:hopfalggen}.

Recall that $\Psi^p=\lambda\Psi$ with $\lambda=\frac{\varepsilon p}{y^{p-1}}$ for some $\varepsilon\in\mathcal{O}_K^{\times}$. Then, $$v_K(\lambda)=e-\frac{p-1}{r}c.$$ 

\subsection{A scaffold on $L/K$}\label{sec:scaffoldp}

The theory of scaffolds was introduced by Byott, Childs and Elder \cite{byottchildselder} as a unification of previous works studying the Hopf-Galois module structure of totally ramified $p$-extensions of local fields with residue characteristic $p$. It is based on the notion of scaffold for such an extension, whose existence gives rise to criteria for the freeness over the associated order in a Hopf-Galois structure in terms of numerical data from the extension. We present here the notion of scaffold and the already mentioned criteria for our extension $L/K$.

Let $b\in\mathbb{Z}$ be coprime with $p$. Call $\mathbb{S}_p=\{0,1,\dots,p-1\}$ and let $\mathfrak{b}\colon\mathbb{S}_{p}\longrightarrow\mathbb{Z}$ be defined by $\mathfrak{b}(s)=bs$. Since $b$ is coprime with $p$, the map $\pi\circ(-\mathfrak{b})$ is bijective, where $\pi\colon\mathbb{Z}\longrightarrow\mathbb{S}_p$ sends each integer number to its remainder mod $p$. Let $\mathfrak{a}\colon\mathbb{S}_p\longrightarrow\mathbb{S}_p$ be its inverse.

\begin{defi}\label{def:scaffold} Let $A$ be a $K$-algebra of dimension $p$ acting $K$-linearly on $L$ and let $\mathfrak{c}\in\mathbb{Z}_{\geq1}$. An $A$-scaffold on $L$ with precision $\mathfrak{c}$ and shift parameter $b$ consists of the following data:
\begin{enumerate}[label=\arabic*)]
    \item\label{def:scaffold1} A family $\{\lambda_n\}_{n\in\mathbb{Z}}$ of elements in $L$ such that $v_L(\lambda_n)=n$ for all $n\in\mathbb{Z}$ and $\lambda_{n_1}\lambda_{n_2}^{-1}\in K$ if and only if $n_1\equiv n_2\,(\mathrm{mod}\,p)$.
    \item\label{def:scaffold2} An element $\Phi\in A$ with $\Phi\cdot1=0$ and the property that for every $n\in\mathbb{Z}$ there is some $u_n\in\mathcal{O}_K^*$ such that:
    $$\begin{cases}
    \Phi\cdot\lambda_n\equiv u_n\lambda_{n+b}\,(\mathrm{mod}\,\lambda_{n+b}\mathfrak{p}_L^{\mathfrak{c}}) &\hbox{ if }\mathfrak{a}(n)\geq1 \\
    \Phi\cdot\lambda_n\equiv0\,(\mathrm{mod}\,\lambda_{n+b}\mathfrak{p}_L^{\mathfrak{c}}) &\hbox{ if }\mathfrak{a}(n)=0 \\
    \end{cases}$$
\end{enumerate}
\end{defi}

We review Elder's result that $L/K$ possesses an $H$-scaffold, where $H$ is the only Hopf-Galois structure on $L/K$. In our case, we take $b$ to be the integer whose existence is due to Theorem \ref{thm:eldertypicalp} applied on $L/K$. Note that since $\mathfrak{a}=(\pi\circ(-\mathfrak{b}))^{-1}$, for each $n\in\mathbb{Z}$, there is $f_n\in\mathbb{Z}$ such that $n=-\mathfrak{a}(n)b+f_np$ (here we identify $\mathfrak{a}(n)=\mathfrak{a}(r(n))$ for $n\in\mathbb{Z}$).

\begin{thm}[\cite{elder2018}, Theorem 3.5]\label{thm:actiononxi} Let $\Psi$ be the $K$-algebra generator of $H$ from Theorem \ref{thm:hopfalggen} and let $$\mathfrak{c}\coloneqq pe-(p-1)\ell.$$ Then, $\Psi\cdot1=0$ and, for a given $1\leq i\leq p-1$, we have $$\Psi\cdot x^i\equiv ix^{i-1}\,(\mathrm{mod}\,x^{i-1}\mathfrak{p}_L^{\mathfrak{c}}).$$
\end{thm}

Note that $0<\ell<\frac{pe}{p-1}$ implies that $\mathfrak{c}>0$. Theorem \ref{thm:actiononxi} gives rise to a scaffold on $L/K$ as follows.

\begin{coro}\label{coro:scaffoldpintr} Let $L/K$ be a typical degree $p$ extension of $p$-adic fields and adopt the notation from Theorem \ref{thm:eldertypicalp}. There is a scaffold on $L/K$ that consists in the following data:
\begin{enumerate}[label=\arabic*)]
    \item The family $\{\lambda_n\}_{n\in\mathbb{Z}}$ of elements in $L$ defined by $\lambda_n=x^{\mathfrak{a}(n)}\pi_K^{f_n}$.
    \item The element $\Psi=-\frac{1}{y}\sum_{m=1}^{p-1}\chi(m)^{-1}\sigma^m\in H$ from Theorem \ref{thm:hopfalggen}.
\end{enumerate}
\end{coro}

The results from \cite{byottchildselder} providing criteria for the freeness are conditioned to the existence of an $A$-scaffold with high enough precision for the intended extension. In our case, we know that $L/K$ possesses an $H$-scaffold, so we can state these criteria without further hypotheses. As a first illustration, \cite[Proposition 2.12]{byottchildselder} becomes as follows.

\begin{pro}\label{pro:hnormbasisscaffold} Let $a'$ be an integer number such that $\mathfrak{a}(a')=p-1$. Then any element $\rho\in L$ with $v_L(\rho)=a'$ is a free generator of $L$ as an $H$-module. Moreover, for each $0\leq h\leq p-1$, $\mathfrak{A}_{L/K}(h)$ is an $\mathcal{O}_K$-order in $H$.
\end{pro}

Note that the last sentence of Proposition \ref{pro:hnormbasisscaffold} is just Proposition \ref{pro:assocordproperties} \ref{pro:assocordproperties1}.

Next, we aim to translate \cite[Theorem 3.1]{byottchildselder} providing criteria for freeness in terms of the $H$-scaffold on $L$ from Corollary \ref{coro:scaffoldpintr}. For each $0\leq h\leq p-1$, denote $\mathbb{S}_p(h)=\{h,h+1,\dots,h+p-1\}$. We seek an integer $a'\in\mathbb{Z}$ such that $$a'\in\mathbb{S}_p(h),\quad\mathfrak{a}(a')=p-1.$$ It is immediate that $$a'=\begin{cases}
    a & \hbox{if }0\leq h\leq a, \\
    a+p & \hbox{if }a+1\leq h\leq p-1
\end{cases}$$ does the job. Now, for each $0\leq i\leq p-1$, let us call \begin{equation}\label{eq:paramscaffold}
    \begin{split}
        d(i;h)&=\Big\lfloor\frac{a'+ib-h}{p}\Big\rfloor,\\w(i;h)&=\mathrm{min}\{d(i+j;h)-d(j;h)\,|\,0\leq j\leq p-1-i\}.
    \end{split}
\end{equation}

In accordance with the notation at \cite[\textsection6]{gildegpgeneral}, let us denote $d(i)=d(i;0)$ and $w(i)=w(i;0)$ for every $0\leq i\leq p-1$. Write $b=pa_0+a$ for the Euclidean division of $b$ by $p$. A straightforward calculation yields $d(p-1)=a+(p-1)a_0$.For each $0\leq i\leq p-1$, let us call $$e_i=\pi_K^{-w(i;h)}\Psi^i.$$

\begin{thm}\label{thm:criteriascaffold} We have that:
\begin{enumerate}[label=\arabic*)]
    \item\label{thm:criteriascaffold1} If $\mathfrak{c}\geq\mathrm{max}(a'-h,1)$, $\{e_i\}_{i=0}^{p-1}$ is an $\mathcal{O}_K$-basis of $\mathfrak{A}_{L/K}(h)$. If in addition $w(i;h)=d(i;h)$ for all $0\leq i\leq p-1$, then $\mathfrak{p}_L^h$ is $\mathfrak{A}_{L/K}(h)$-free.
    \item\label{thm:criteriascaffold2} If $\mathfrak{c}\geq p+a'-h$, $\mathfrak{p}_L^h$ is $\mathfrak{A}_{L/K}(h)$-free if and only if $w(i;h)=d(i;h)$ for all $0\leq i\leq p-1$.
\end{enumerate}
\end{thm}

The inequality at \ref{thm:criteriascaffold1} is called the weak condition, while the one at \ref{thm:criteriascaffold2} is the strong one. In order to explore whether they are fulfilled, we work with an alternative form of the precision $\mathfrak{c}$.

\begin{pro}[\cite{gildegpgeneral}, Proposition 6.1]\label{pro:precscaffold} Let $L/K$ be a typical degree $p$ extension of $p$-adic fields. The precision of the scaffold on $L/K$ can be written as $$\mathfrak{c}=p\Big(v_K(\lambda)-d(p-1)\Big)+a.$$ In particular, the remainder of $\mathfrak{c}$ mod $p$ is $a$.
\end{pro}

One of the key facts that allow to apply the theory of scaffolds in the case $h=0$ is that $\mathfrak{c}\geq p+a$ if and only if $\ell<\frac{pe}{p-1}-1$ (see \cite[Proposition 6.3]{gildegpgeneral}). In the following, we expand the relationship between a lower bound for $\mathfrak{c}$ and an upper bound for $\ell$.

\begin{lem}\label{lem:relprecramjump} Let $k\in\mathbb{Z}$ with $0\leq k<\frac{pe}{p-1}$. If $\mathfrak{c}\geq kp+a$, then $\ell<\frac{pe}{p-1}-k$. Moreover, if $k\leq p-1$ and $a<p-k+1$, then the converse holds.
\end{lem}
\begin{proof}
Suppose that $\mathfrak{c}\geq kp+a$. By definition, $\mathfrak{c}=pe-(p-1)\ell$, so $$\ell=\frac{pe}{p-1}-\frac{\mathfrak{c}}{p-1}\leq\frac{pe}{p-1}-\frac{kp+a}{p-1}<\frac{pe}{p-1}-k.$$ 

Next, assume that $\ell<\frac{pe}{p-1}-k$, $k\leq p-1$, $a<p-k+1$ and $\ell<\frac{pe}{p-1}-k$. From the latter we obtain that $\mathfrak{c}>k(p-1)=kp-k$. Since $k\leq p-1$, $\mathfrak{c}>(k-1)p$. Suppose that $\mathfrak{c}<kp$ (which is possible only for $k\geq2$). By Proposition \ref{pro:precscaffold}, $\mathfrak{c}\equiv a(\mathrm{mod}\,p)$, so $\mathfrak{c}=(k-1)p+a$. Now, the inequality $a<p-k+1$ yields $\mathfrak{c}=(k-1)p+a<kp-k+1=k(p-1)+1$, that is, $\mathfrak{c}\leq k(p-1)$, which is a contradiction. We get $\mathfrak{c}\geq kp$, whence $\mathfrak{c}\geq kp+a$.
\end{proof}

\begin{rmk}\normalfont The conditions that $k\leq p-1$ and $a<p-k+1$ are automatically satisfied for $k\in\{0,1\}$. Thus, for these values of $k$, the converse holds unconditionally. For $k=0$, we obtain a tautology (both $\mathfrak{c}\geq a$ and $\ell<\frac{pe}{p-1}$ hold unconditionally). As for $k=1$, we recover \cite[Proposition 6.3]{gildegpgeneral}. On the other hand, for $k=2$, Lemma \ref{lem:relprecramjump} says that $\mathfrak{c}\geq2p+a$, then $\ell<\frac{pe}{p-1}-2$, the converse holding if $a<p-1$.
\end{rmk}

\begin{coro}\label{coro:criteriamidtransl} Given $0\leq h\leq p-1$, the following statements are satisfied:
\begin{enumerate}[label=\arabic*)]
    \item\label{coro:criteriamidtransl1} If $0\leq h\leq a$ or $\ell<\frac{pe}{p-1}-1$, then $\mathfrak{A}_{L/K}(h)$ has $\mathcal{O}_K$-basis $\{e_i\}_{i=0}^{p-1}$ and $w(i;h)=d(i;h)$ for all $0\leq i\leq p-1$ implies that $\mathfrak{p}_L^h$ is $\mathfrak{A}_{L/K}(h)$-free.
    \item\label{coro:criteriamidtransl2} If $0\leq h\leq a$ and $\ell<\frac{pe}{p-1}-1$ or $a+1\leq h\leq p-1$ and $\ell<\frac{pe}{p-1}-2$, then $\mathfrak{p}_L^h$ is $\mathfrak{A}_{L/K}(h)$-free if and only if $w(i;h)=d(i;h)$ for all $0\leq i\leq p-1$.
\end{enumerate}
\end{coro}
\begin{proof}
Suppose that $0\leq h\leq a$, so $a'=a$. By Proposition \ref{pro:precscaffold}, $\mathfrak{c}=p(v_K(\lambda)-d(p-1))+a\geq a\geq\mathrm{max}(a-h,1)$, so the weak condition is satisfied. We obtain that $\mathfrak{A}_{L/K}(h)$ has $\mathcal{O}_K$-basis $\{\pi_K^{-w(i;h)}\Psi^i\}_{i=0}^{p-1}$ and $w(i;h)=d(i;h)$ for all $0\leq i\leq p-1$ implies that $\mathfrak{p}_L^h$ is $\mathfrak{A}_{L/K}(h)$-free. In addition, the strong condition $\mathfrak{c}\geq p+a-h$ holds if and only if $\mathfrak{c}\geq p+a$. By Lemma \ref{lem:relprecramjump}, this is equivalent to $\ell<\frac{pe}{p-1}-1$. Under this hypothesis, $w(i;h)=d(i;h)$ for all $0\leq i\leq p-1$ if and only if $\mathfrak{p}_L^h$ is $\mathfrak{A}_{L/K}(h)$-free.

Now assume that $a+1\leq h\leq p-1$, in which case $a'=a+p$. The weak condition is now $\mathfrak{c}\geq p+a-h$. Once again, this holds if and only if $\ell<\frac{pe}{p-1}-1$. On the other hand, the strong condition becomes $\mathfrak{c}\geq 2p+a-h$, and this is satisfied if and only if $\mathfrak{c}\geq 2p+a$. Assume that $\ell<\frac{pe}{p-1}-2$. Note that $a<p-1$ because $a+1\leq p-1$. Then we can apply Lemma \ref{lem:relprecramjump} with $k=2$ to obtain that $\mathfrak{c}\geq2p+a$. Then the strong condition is indeed satisfied, and $\mathfrak{p}_L^h$ is $\mathfrak{A}_{L/K}(h)$-free if and only if $w(i;h)=d(i;h)$ for all $0\leq i\leq p-1$.
\end{proof}

In particular, Corollary \ref{coro:criteriamidtransl} says that the weak condition always holds if $0\leq h\leq a$ and it does under the assumption that $\ell<\frac{pe}{p-1}-1$ if $a+1\leq h\leq p-1$. The only case that is not covered by Corollary \ref{coro:criteriamidtransl} corresponds to $a+1\leq h\leq p-1$ and $\ell\geq\frac{pe}{p-1}-1$, which will be treated in Section \ref{sec:exceptionalcase}. For now, we shall assume that the weak condition holds. Note that it can be restated as $0\leq h\leq a$ or $\ell<\frac{pe}{p-1}-1$.

Under our assumption, \ref{coro:criteriamidtransl1} gives that $\mathfrak{A}_{L/K}(h)$ has $\mathcal{O}_K$-basis $\{e_i\}_{i=0}^{p-1}$, where $e_i=\pi_K^{-w(i;h)}\Psi^i$. The idea is that the negatives of the numbers $w(i;h)$ are the powers by which we normalize the elements $\Psi^i$ so as to obtain an $\mathcal{O}_K$-basis of $\mathfrak{A}_{L/K}(h)$. 

We can interpret the numbers $d(i;h)$ in an analogous way. Call $\theta=\pi_L^{a'}$ and $\mathfrak{A}_{\theta}(h)=\{\gamma\in H\,\mid\,\gamma\cdot\theta\in\mathfrak{p}_L^h\}$. For each $0\leq i\leq p-1$, call $$f_i\coloneqq\pi_K^{-d(i;h)}\Psi^i.$$

\begin{pro}\label{pro:basisathetah} For every $0\leq h\leq p-1$, $\{f_i\}_{i=0}^{p-1}$ is an $\mathcal{O}_K$-basis of $\mathfrak{A}_{\theta}(h)$.
\end{pro}
\begin{proof}
Let $\gamma\in H$ and write $\gamma=\sum_{i=0}^{p-1}\gamma_i\Psi^i$ with $\gamma_i\in K$. Then $\gamma\in\mathfrak{A}_{\theta}(h)$ if and only if $\gamma\cdot\theta\in\mathfrak{p}_L^h$; if and only if $\sum_{i=0}^{p-1}\gamma_i\Psi^i\cdot\theta\in\mathfrak{p}_L^h$. Now, $$v_L(\gamma_i\Psi^i\cdot\theta)=pv_K(\gamma_i)+a'+ib\equiv(i+1)a\,(\mathrm{mod}\,p).$$ Hence, the $L$-valuations of the elements $\gamma_i\Psi^i\cdot\theta$ are distinct modulo $p$, so $\sum_{i=0}^{p-1}\gamma_i\Psi^i\cdot\theta\in\mathfrak{p}_L^h$ if and only if $\gamma_i\Psi^i\cdot\theta\in\mathfrak{p}_L^h$ for all $0\leq i\leq p-1$. Now, for each $0\leq i\leq p-1$, $\gamma_i\Psi^i\cdot\theta\in\mathfrak{p}_L^h$ if and only if $$pv_K(\gamma_i)+a'+ib\geq h,$$ which is equivalent to $$v_K(\gamma_i)\geq-\frac{a'+ib-h}{p},$$ and this clearly happens if and only if $v_K(\gamma_i)\geq-d(i;h)$. We have proved that $\gamma\in\mathfrak{A}_{\theta}(h)$ if and only if $\gamma_i\in\pi_K^{-d(i;h)}\mathcal{O}_K$ for all $0\leq i\leq p-1$, so the statement follows.
\end{proof}

\section{The equalities $d(i;h)=w(i;h)$}\label{sec:equaldiwi}

Let $L/K$ be a typical degree $p$ extension of $p$-adic fields. Suppose that $0\leq h\leq a$ or $\ell<\frac{pe}{p-1}-1$, so that we are in the situation of Corollary \ref{coro:criteriamidtransl} \ref{coro:criteriamidtransl1}. Then, we know that $\mathfrak{A}_{L/K}(h)$ has $\mathcal{O}_K$-basis $\{e_i\}_{i=0}^{p-1}$. Moreover, the equalities $d(i;h)=w(i;h)$ for $0\leq i\leq p-1$ are a sufficient condition for the $\mathfrak{A}_{L/K}(h)$-freeness of $\mathfrak{p}_L^h$, that become necessary if the strong condition of the scaffold holds. If $h=0$, the strong condition does not hold when $\ell\geq\frac{pe}{p-1}-1$ (that is, $L/K$ is almost maximally ramified), and it is seen that the validity of the equalities $d(i;h)=w(i;h)$ is not a necessary condition. Instead, $\mathcal{O}_L$ is $\mathfrak{A}_{L/K}$-free if and only if the length of the continued fraction expansion of $\frac{b}{p}$ is at most $4$, which is a less restrictive condition (see \cite[\textsection 7]{gildegpgeneral}). Nevertheless, if $1\leq h\leq p-1$, we will see that the validity of these equalities is still a necessary condition, and therefore a characterization, for the $\mathfrak{A}_{L/K}(h)$-freeness of $\mathfrak{p}_L^h$. In other words, we shall prove the following result.

\begin{thm}\label{thm:charactfrenessnumeric} Suppose that $h>0$ and that $1\leq h\leq a$ or $\ell<\frac{pe}{p-1}-1$. Then, $\mathfrak{p}_L^h$ is $\mathfrak{A}_{L/K}(h)$-free if and only if $$w(i;h)=d(i;h)\quad\hbox{for all }0\leq i\leq p-1.$$
\end{thm}

We already know the validity of this result when the strong condition holds; namely, $\ell<\frac{pe}{p-1}-1$ if $1\leq h\leq a$, and $\ell<\frac{pe}{p-1}-2$ if $a+1\leq h\leq p-1$. Moreover, under the weak condition, we know the implication from right to left. Thus, we shall assume that the weak condition holds but the strong does not, and prove the implication from left to right. Namely: \begin{equation}\label{eq:hypothesisweaknotstr}
    \begin{cases}
        \frac{pe}{p-1}-1\leq\ell<\frac{pe}{p-1} & \hbox{if }1\leq h\leq a, \\
        \frac{pe}{p-1}-2\leq\ell<\frac{pe}{p-1}-1 & \hbox{if }a+1\leq h\leq p-1.
    \end{cases}
\end{equation}

\begin{lem}\label{lem:diffequal1} Suppose that $1\leq h\leq p-1$ and that \eqref{eq:hypothesisweaknotstr} holds. Then, $$v_K(\lambda)-w(p-1;h)=1.$$
\end{lem}
\begin{proof}
Note that, under \eqref{eq:hypothesisweaknotstr}, we have from Lemma \ref{lem:relprecramjump} that $\mathfrak{c}=a'$. Using Proposition \ref{pro:precscaffold}, $$v_K(\lambda)-d(p-1)=\begin{cases}
    0 & \hbox{if }1\leq h\leq a, \\
    1 & \hbox{if }a+1\leq h\leq p-1.
\end{cases}$$ On the other hand, by definition, $w(p-1;h)=d(p-1;h)-d(0;h)$. Since $a'\in\mathbb{S}_p(h)$, $$d(0;h)=\Big\lfloor\frac{a'-h}{p}\Big\rfloor=0.$$ Thus, $w(p-1;h)=d(p-1;h)$.

If $1\leq h\leq a$, we have \begin{equation*}
    \begin{split}
        d(p-1;h)&=\Big\lfloor\frac{a+(p-1)b-h}{p}\Big\rfloor
        \\&=(p-1)a_0+a+\Big\lfloor\frac{-h}{p}\Big\rfloor
        \\&=d(p-1)-1,
    \end{split}
\end{equation*} the last equality due to $h>0$. Hence, $$v_K(\lambda)-w(p-1;h)=v_K(\lambda)-d(p-1)+1=1.$$ Otherwise, if $a+1\leq h\leq p-1$, $$d(p-1;h)=\Big\lfloor\frac{a+(p-1)b-h}{p}\Big\rfloor+1=d(p-1).$$ Then, $$v_K(\lambda)-w(p-1;h)=v_K(\lambda)-d(p-1)=1,$$ as we wanted.
\end{proof}

Let us call $A_h=\mathfrak{A}_{L/K}(h)$. Under the weak condition (which is implied by \eqref{eq:hypothesisweaknotstr}), we know that $A_h$ has $\mathcal{O}_K$-basis $\{e_i\}_{i=0}^{p-1}$. In addition, recall that $A_h$ is an $\mathcal{O}_K$-order; in particular a ring. It is straightforward to check that $\pi_KA_h$ is an ideal of $A_h$. Consider the quotient ring $\overline{A_h}\coloneqq A_h/\pi_KA_h$, which is naturally endowed with $\overline{K}$-vector space structure. Inside this quotient, define $J\coloneqq\langle\{\overline{e_i}\}_{i=1}^{p-1}\rangle_{\overline{K}}$. Since $\{\overline{e_i}\}_{i=1}^{p-1}$ is a $\overline{K}$-basis of $\overline{A_h}$ and $\Psi^p=\lambda\Psi$ from Theorem \ref{thm:hopfalggen}, $e_ie_j\in J$ for every $0\leq i,j\leq p-1$ such that $i>0$ or $j>0$. We find that $J$ is an ideal of $\overline{A_h}$. Write $\mathrm{Jac}(\overline{A_h})$ for the Jacobson radical of $\overline{A_h}$.

\begin{lem}\label{lem:jacobsonAbar} Suppose that $1\leq h\leq p-1$ and that \eqref{eq:hypothesisweaknotstr} holds. Then $$J=\mathrm{Jac}(\overline{A_h}).$$
\end{lem}
\begin{proof}
First, we prove that $J$ is a nilpotent ideal of $\overline{A_h}$, which implies that $J\subseteq\mathrm{Jac}(\overline{A_h})$. Given $0\leq i\leq p-1$, one has $e_i^{p-1}=\pi_K^{-(p-1)w(i;h)}\Psi^{i(p-1)}$. Now, $$\Psi^{i(p-1)}=\Psi^{(i-1)p+p-i}=(\Psi^p)^{i-1}\Psi^{p-i}=\lambda^{i-1}\Psi^{p-1}.$$ Therefore, \begin{equation*}
    \begin{split}
        e_i^{p-1}&=\lambda^{i-1}\pi_K^{-(p-1)w(i;h)}\Psi^{p-1}
        \\&=\lambda^{i-1}\pi_K^{w(p-1;h)-(p-1)w(i;h)}e_{p-1}.
    \end{split}
\end{equation*} Since $A_h$ is a ring, the coefficient has non-negative valuation. That is, $$(i-1)v_K(\lambda)+w(p-1;h)-(p-1)w(i;h)\geq0.$$ On the other hand, \begin{equation*}
    \begin{split}
        e_i^p&=\pi_K^{-pw(i;h)}\Psi^{ip}
        \\&=\lambda^i\pi_K^{-pw(i;h)}\Psi^i
        \\&=\lambda^i\pi_K^{-(p-1)w(i;h)}e_i.
    \end{split}
\end{equation*} The $K$-valuation of the coefficient in this relation is \begin{equation*}
    \begin{split}
        iv_K(\lambda)-(p-1)w(i;h)&=[(i-1)v_K(\lambda)+w(p-1;h)-(p-1)w(i;h)]+[v_K(\lambda)-w(p-1;h)]
        \\&\geq v_K(\lambda)-w(p-1;h)=1,
    \end{split}
\end{equation*} where we have used Lemma \ref{lem:diffequal1}. Thus, passage to the quotient $A_h/\pi_KA_h$ annihilates the coefficient. We deduce that $\overline{e_i}^p=\overline{0}$ for all $1\leq i\leq p-1$. Since $A_h$ is commutative, we deduce that $J$ is nilpotent.

On the other hand, since $\overline{A_h}$ is the $\overline{K}$-span of $\overline{e_0},\dots,\overline{e_{p-1}}$ and $J$ is the $\overline{K}$-span of $\overline{e_1},\dots,\overline{e_{p-1}}$, there are ring isomorphisms $\overline{A_h}/J\cong\overline{K}e_0\cong\overline{K}$. Therefore, $\overline{A_h}/J$ is a field, so $J$ is a maximal ideal of $\overline{A_h}$, whence $\mathrm{Jac}(\overline{A_h})\subseteq J$. This concludes the proof.
\end{proof}

\begin{pro}\label{pro:hardimplication} Suppose that $1\leq h\leq p-1$ and that \eqref{eq:hypothesisweaknotstr} holds. If $\mathfrak{p}_L^h$ is $\mathfrak{A}_{L/K}(h)$-free, then $w(i;h)=d(i;h)$ for all $0\leq i\leq p-1$.
\end{pro}
\begin{proof}
Let $\theta=\pi_L^{a'}$ and call $I_h=\mathfrak{A}_{\theta}(h)$, $\overline{I_h}=I_h/\pi_KI_h$. We already know that $I_h$ is an $A_h$-fractional ideal. Moreover, by Proposition \ref{pro:hnormalbasisfract}, the hypothesis that $\mathfrak{p}_L^h$ is $A_h$-free is equivalent to $I_h$ being principal as an $A_h$-fractional ideal. Let $\alpha\in I_h$ be an $A_h$-generator of $I_h$ and let $\psi_{\alpha}\colon A_h\longrightarrow I_h$ be the map defined by $\psi_{\alpha}(\lambda)=\lambda\alpha$. Since $I_h$ is $A_h$-principal, $\psi_{\alpha}$ is an $A_h$-linear isomorphism. Moreover, $\psi_{\alpha}(\pi_KA_h)=\pi_KI_h$. Therefore, there is an $\overline{A_h}$-linear isomorphism $\overline{A_h}\cong\overline{I_h}$. Equivalently, $\overline{I_h}$ is free as an $\overline{A_h}$-module. Under this isomorphism, the ideal $J$ of $\overline{A_h}$ is sent to $J\overline{I_h}$. Hence, $\overline{I_h}/J\overline{I_h}\cong\overline{A_h}/J\cong\overline{K}$ as $\overline{K}$-vector spaces, so $\mathrm{dim}_{\overline{K}}(\overline{I_h}/J\overline{I_h})=1$. In other words, any non-zero element of $\overline{I_h}/J\overline{I_h}$ generates it as a $\overline{K}$-vector space. 

Recall that $J$ is generated by $\overline{e_1},\dots,\overline{e_{p-1}}$, while $\overline{I_h}$ is generated by $\overline{f_0},\dots,\overline{f_{p-1}}$ by Proposition \ref{pro:basisathetah}. We claim that $\overline{f_0}\notin J\overline{I_h}$. Indeed, for $1\leq i\leq p-1$ and $0\leq j\leq p-1$, $\overline{e_i}\overline{f_j}=\overline{\pi_K^{-w(i;h)-d(j;h)}\Psi^{i+j}}$ is a linear combination of $\overline{f_0},\dots,\overline{f_{p-1}}$. Since $\overline{f_0}=1_{\overline{A}}$ and $\Psi^p=\lambda\Psi$, $\overline{f_0}$ cannot arise as a linear combination of elements of the form $\overline{e_i}\overline{f_j}$.

Thus, we have that $\overline{I_h}/J\overline{I_h}$ is generated by $\overline{f_0}$ as a $\overline{K}$-vector space. Since $\overline{A_h}/J\cong\overline{K}$, we get $\overline{I_h}=\overline{A_h}\overline{f_0}+J\overline{I_h}$. By Lemma \ref{lem:jacobsonAbar}, $J=\mathrm{Jac}(\overline{A_h})$. Applying Nakayama's lemma, we obtain that $\overline{I_h}=\overline{A_h}\overline{f_0}$, that is, $\overline{f_0}$ is a free generator of $\overline{I_h}$ as an $\overline{A_h}$-module. Now, given $0\leq i\leq p-1$, $$e_if_0=e_i=\pi_K^{-w(i;h)}\Psi^i=\pi_K^{d(i;h)-w(i;h)}f_i,$$ and passage to the quotient gives $$\overline{e_if_0}=\overline{\pi_K^{d(i;h)-w(i;h)}f_i}.$$ Since $\overline{I_h}=\overline{A_h}\overline{f_0}$, for each $0\leq i\leq p-1$ we must have $\overline{e_if_0}\neq\overline{0}$, so $w(i;h)=d(i;h)$. This concludes the proof.
\end{proof}

Theorem \ref{thm:charactfrenessnumeric} follows immediately from Corollary \ref{coro:criteriamidtransl} and Proposition \ref{pro:hardimplication}.

\section{The equality $\mathfrak{A}_{L/K}(h)=\mathfrak{A}_{\theta}(h)$}\label{sec:equalassocorders}

We continue with the same setting as in the previous section. For each $1\leq h\leq p-1$, we now know that $\mathfrak{p}_L^h$ is $\mathfrak{A}_{L/K}(h)$-free if and only if $w(i;h)=d(i;h)$ for all $0\leq i\leq p-1$. On the other hand, for $h=0$, it is so unless $\ell\geq\frac{pe}{p-1}-1$. Moreover, in that case, $w(i)=d(i)$ for all $0\leq i\leq p-1$ if and only if $a\mid p-1$ (see \cite[Proposition 6.11]{gildegpgeneral}). In this section we are going to characterize the validity of the equalities $w(i;h)=d(i;h)$, $0\leq i\leq p-1$, for arbitrary $0\leq h\leq p-1$. Unlike in the previous section, we will find uniform behavior for all such $h$. 

Note that the elements $e_i=\pi_K^{-w(i;h)}\Psi^i$ form an $\mathcal{O}_K$-basis of $\mathfrak{A}_{L/K}(h)$ by Corollary \ref{coro:criteriamidtransl} \ref{coro:criteriamidtransl1}, while the elements $f_i=\pi_K^{-d(i;h)}\Psi^i$ form an $\mathcal{O}_K$-basis of $\mathfrak{A}_{\theta}(h)$ for $\theta=\pi_L^{a'}$. Therefore, the validity of the already mentioned equalities is equivalent to $\mathfrak{A}_{L/K}(h)=\mathfrak{A}_{\theta}(h)$. Using Proposition \ref{pro:hnormalbasisfract} \ref{pro:hnormalbasisfract3}, this equality holds if and only if $\mathfrak{A}_{\theta}(h)$ is a ring. This is the condition we shall work with. First, we can characterize it in terms of the numbers $d(i;h)$ by making use of Proposition \ref{pro:basisathetah}.

\begin{coro}\label{coro:twoimpl} Given $0\leq h\leq p-1$, $\mathfrak{A}_{\theta}(h)$ is a ring if and only if for every $0\leq i,j\leq p-1$ the following two implications hold:
\begin{enumerate}[label=\arabic*)]
    \item\label{coro:twoimpl1} If $i+j\leq p-1$, then $d(i;h)+d(j;h)\leq d(i+j;h)$.
    \item\label{coro:twoimpl2} If $i+j\geq p$, then $d(i;h)+d(j;h)\leq v_K(\lambda)+d(i+j+1-p;h)$.
\end{enumerate}
\end{coro}
\begin{proof}
Let $0\leq i,j\leq p-1$ such that $i+j\leq p-1$. It is trivial that the product $\pi_K^{-d(i;h)}\Psi^i\pi_K^{-d(j;h)}\Psi^j$ belongs to $\mathfrak{A}_{\theta}(h)$ if and only if $d(i;h)+d(j;h)\leq d(i+j;h)$. Now, assume that $i+j\geq p$. We have that \begin{equation*}
    \begin{split}
        \pi_K^{-d(i;h)}\Psi^i\pi_K^{-d(j;h)}\Psi^j&=\pi_K^{-d(i;h)-d(j;h)}\Psi^{i+j}
        \\&=\lambda\pi_K^{-d(i;h)-d(j;h)}\Psi^{i+j+1-p}.
    \end{split}
\end{equation*} This belongs to $\mathfrak{A}_{\theta}(h)$ if and only if it sends $\theta$ to $\mathfrak{p}_L^h$. Taking $L$-valuations, we have $$v_L(\lambda)-p(d(i;h)+d(j;h))+a'+(i+j+1-p)b\geq h,$$ that is, $$p(d(i;h)+d(j;h))\leq v_L(\lambda)+a'+(i+j+1-p)b-h.$$ This is clearly equivalent to the inequality at \ref{coro:twoimpl2}.
\end{proof}

We shall study whether the implications at Corollary \ref{coro:twoimpl} hold so as to characterize the equality $\mathfrak{A}_{L/K}(h)=\mathfrak{A}_{\theta}(h)$.

\subsection*{Interlude: fractional parts and continued fractions}

Before continuing, let us introduce some notation. For a rational number $\alpha\in\mathbb{Q}$, we write $\widehat{\alpha}$ for its fractional part, so that $\alpha=\lfloor \alpha\rfloor+\widehat{\alpha}$. Also, we write $||x||$ for the distance to the nearest integer of $x$. If $m$ is a positive integer and $x\in\mathbb{Z}$, we write $[x]_m$ for the remainder of $x$ mod $m$. On this way, $x=\lfloor\frac{x}{m}\rfloor m+[x]_m$. If $p\nmid x$, we write $x^{-1}$ for the remainder mod $p$ of the modular inverse of $x$ in $\mathbb{Z}/p\mathbb{Z}$.

Recall that we write $$\frac{b}{p}=[a_0;a_1,\dots,a_n]$$ for the continued fraction expansion of $\frac{b}{p}$. Note that the continued fraction expansion of $\frac{a}{p}$ is $$\frac{a}{p}=[0;a_1,\dots,a_n].$$ We shall write $$\frac{p_i}{q_i}=[0;a_1,\dots,a_i],\quad1\leq i\leq n,$$ for the convergents of the continued fraction expansion of $\frac{a}{p}$. We are interested in the denominators $\{q_i\}_{i=1}^n$; consequently we will refer to them as the convergents. Note that they are the same in the continued fraction expansions of $\frac{a}{p}$ and $\frac{b}{p}$, so we can think of the $\{q_i\}_{i=0}^n$ indistinctly as the convergents of one or the other. It is easy to see that $q_0=1$, $q_1=a_1$ and $q_n=p$. Moreover, $$q_i=a_iq_{i-1}+q_{i-2},\quad 2\leq i\leq n.$$ We also recall some facts that will be used in the sequel (for the proofs, see \cite[I, \textsection2]{lang1995}).

\begin{pro}\label{pro:contfrac} Let $a\in\mathbb{Z}$ with $0<a<p$ and write $\frac{a}{p}=[a_0;a_1,\dots,a_n]$ for the continued fraction expansion of $\frac{a}{p}$. Let $\frac{p_i}{q_i}$ be its $i$-th convergent, $0\leq i\leq n$.
\begin{enumerate}[label=\arabic*)]
    \item\label{pro:contfrac1} $||q_{i+1}\frac{a}{p}||<||q_i\frac{a}{p}||$ for every $1\leq i\leq n-1$.
    \item\label{pro:contfrac2} If $q\in\mathbb{Z}$ and $0<q<q_i$, then $||q\frac{a}{p}||\geq||q_{i-1}\frac{a}{p}||$, and in particular, $||q\frac{a}{p}||>||q_i\frac{a}{p}||$.
    \item\label{pro:contfrac3} If $1\leq i\leq n-1$ is even (resp. odd), then $\widehat{q_i\frac{a}{p}}<\frac{1}{2}$ (resp. $\widehat{q_i\frac{a}{p}}>\frac{1}{2}$).
    \item\label{pro:contfrac4} $||q_{n-1}\frac{a}{p}||=\frac{1}{p}$, so $\widehat{q_{n-1}\frac{a}{p}}=\frac{1}{p}$ if $n$ is odd and $\widehat{q_{n-1}\frac{a}{p}}=\frac{p-1}{p}$ if $n$ is even.
\end{enumerate}
\end{pro}

\subsection{The case $0\leq h\leq a$}

Let us characterize the equality $\mathfrak{A}_{\theta}(h)=\mathfrak{A}_{L/K}(h)$ when $0\leq h\leq a$. To this end, we study the validity (or otherwise) of the two implications at Corollary \ref{coro:twoimpl}. 

Recall that in the case $0\leq h\leq a$ we have that $a'=a$. Hence, for each $0\leq k\leq p-1$, we have $$d(k;h)=\Big\lfloor\frac{a+kb-h}{p}\Big\rfloor=ka_0+\Big\lfloor\frac{(k+1)a-h}{p}\Big\rfloor.$$ For every such a $k$, call $$r_k\coloneqq[(k+1)a-h]_p.$$

\begin{pro}\label{pro:ineqremainders} Suppose that $0\leq h\leq a$. Then, given $0\leq i,j\leq p-1$, the following statements hold:
\begin{enumerate}[label=\arabic*)]
    \item\label{pro:ineqremainders1} If $i+j\leq p-1$, then $d(i;h)+d(j;h)\leq d(i+j;h)$ if and only if \begin{equation}\label{eq:ineqremainders}
    r_i+r_j\geq a-h.
\end{equation}
    \item\label{pro:ineqremainders2} If $i+j\geq p$, then $d(i;h)+d(j;h)\leq v_K(\lambda)+d(i+j+1-p;h)$.
\end{enumerate}
\end{pro}
\begin{proof}
Fix $0\leq i,j\leq p-1$.
\begin{enumerate}[label=\arabic*)]
    \item Rewriting integer parts in terms of fractional parts, we see that $d(i;h)+d(j;h)\leq d(i+j;h)$ if and only if the inequality \begin{equation}\label{eq:ineqfractpart}
    \widehat{\frac{(i+1)a-h}{p}}+\widehat{\frac{(j+1)a-h}{p}}\geq\frac{a-h}{p}+\widehat{\frac{(i+j+1)a-h}{p}}
\end{equation} is satisfied. Let us note that $\widehat{\frac{x}{p}}=\frac{[x]_p}{p}$ for every $0\leq x\leq p-1$. Thus, multiplying by $p$, we see that \eqref{eq:ineqfractpart} is equivalent to $$r_i+r_j\geq a-h+r_{i+j}.$$ Now, note that $r_{i+j}=[r_i+r_j-(a-h)]_p$. Indeed, \begin{equation*}
    \begin{split}
        [r_i+r_j-(a-h)]_p&=[ia+a-h+ja+a-h-(a-h)]_p
        \\&=[(i+j)a+a-h]_p
        \\&=r_{i+j}.
    \end{split}
\end{equation*} If we carry out the Euclidean division of $r_i+r_j-(a-h)$ by $p$, we get $$r_i+r_j=\Big\lfloor\frac{r_i+r_j-(a-h)}{p}\Big\rfloor p+a-h+r_{i+j}.$$ Thus, $r_i+r_j\geq a-h+r_{i+j}$ if and only if $\Big\lfloor\frac{r_i+r_j-(a-h)}{p}\Big\rfloor$ is non-negative. This is clearly equivalent to $r_i+r_j\geq a-h$, as we wanted.
    \item If $i+j\geq p$, we have \begin{equation*}
    \begin{split}
        d(i;h)+d(j;h)&=(i+j)a_0+\Big\lfloor\frac{(i+1)a-h}{p}\Big\rfloor+\Big\lfloor\frac{(j+1)a-h}{p}\Big\rfloor
        \\&\leq(i+j)a_0+\Big\lfloor\frac{(i+j+2)a-h}{p}\Big\rfloor
        \\&=(i+j)a_0+a+\Big\lfloor\frac{(i+j+2-p)a-h}{p}\Big\rfloor
        \\&=a+(p-1)a_0+(i+j+1-p)a_0+\Big\lfloor\frac{(i+j+2-p)a-h}{p}\Big\rfloor
        \\&=d(p-1)+\Big\lfloor\frac{a+(i+j+1-p)b-h}{p}\Big\rfloor
        \\&\leq v_K(\lambda)+d(i+j+1-p;h).
    \end{split}
\end{equation*}
\end{enumerate}
\end{proof}

Thus, for $0\leq h\leq a$, $\mathfrak{A}_{L/K}(h)=\mathfrak{A}_{\theta}(h)$ is equivalent to the validity of the inequality \eqref{eq:ineqremainders}. Next, we characterize this inequality in terms of $h$, so as to obtain the desired result.

\begin{pro}\label{pro:characteqathetah} Suppose that $0\leq h\leq a$. Then, the following statements hold:
\begin{enumerate}[label=\arabic*)]
    \item\label{pro:characteqathetah1} If $a=1$, then $\mathfrak{A}_{\theta}(h)=\mathfrak{A}_{L/K}(h)$ for all $h\in\{0,1\}$.
    \item\label{pro:characteqathetah2} If $a>1$, then $\mathfrak{A}_{\theta}(h)=\mathfrak{A}_{L/K}(h)$ if and only if:
    \begin{enumerate}[label=\arabic*)]
        \item $h\in\{a-a_n,a\}$, if $n$ is even.
        \item $a-\frac{a_n}{2}\leq h\leq a$, if $n$ is odd.
    \end{enumerate}
\end{enumerate}
\end{pro}
\begin{proof}
By Corollary \ref{coro:twoimpl} and Proposition \ref{eq:ineqremainders}, $\mathfrak{A}_{\theta}(h)=\mathfrak{A}_{L/K}(h)$ if and only if \eqref{eq:ineqremainders} holds. We will test the validity (or otherwise) of the intended equality through that condition.

First, note that \eqref{eq:ineqremainders} holds trivially if $h=a$. For $a=1$, in order to prove \ref{pro:characteqathetah1}, it remains to validate the equality $\mathfrak{A}_{\theta}(h)=\mathfrak{A}_{L/K}(h)$ for $h=0$. Note that this corresponds to the equality $\mathfrak{A}_{L/K}=\mathfrak{A}_{\theta}$, because $\mathfrak{p}_L^h=\mathcal{O}_L$, and, since $a=1$ divides $p-1$, we can invoke \cite[Proposition 6.11]{gildegpgeneral}. Alternatively, if we replace $a=1$ and $h=0$ at \eqref{eq:ineqremainders}, we get $$[i+1]_p+[j+1]_p\geq1.$$ This is true unless $[i+1]_p=[j+1]_p=0$; but in that case we have that $i=j=p-1$, contradicting that $i+j\leq p-1$.

In what remains, we assume that $a>1$. Necessarily, $n>1$. Set $q\coloneqq[a^{-1}h-1]_p$. Equivalently, $q$ is the integer number between $0$ and $p-1$ defined by $(q+1)a\equiv h\,(\mathrm{mod}\,p)$. 

Suppose first that $n$ is even. By Proposition \ref{pro:contfrac} \ref{pro:contfrac4}, $[q_{n-1}a]_p=p-1$, so $q_{n-1}a\equiv-1\,(\mathrm{mod}\,p)$. Note that since $p=q_n=a_nq_{n-1}+q_{n-2}$, we get that $q_{n-2}a\equiv a_n\,(\mathrm{mod}\,p)$.


Let us assume that $h\notin\{a-a_n,a\}$. If $q+q_{n-1}<p$, for $i=q$ and $j=q_{n-1}$ we have $i+j\leq p-1$ and $$r_i+r_j=[(q+1)a-h]_p+[(q_{n-1}+1)a-h]_p=a-h-1<a-h,$$ so \eqref{eq:ineqremainders} is not satisfied. On the other hand, we claim that the case $q+q_{n-1}\geq p$ is not possible if $h>a-a_n$. From $(q+1)a\equiv h\,(\mathrm{mod}\,p)$, we get $qa\equiv-(a-h)\,(\mathrm{mod}\,p)$. Multiplying by $q_{n-1}$, we get $q\equiv(a-h)q_{n-1}\,(\mathrm{mod}\,p)$. Now, $0<(a-h)q_{n-1}<a_nq_{n-1}<p$, whence $q=(a-h)q_{n-1}$. Therefore, $q+q_{n-1}=(a-h+1)q_{n-1}\leq a_nq_{n-1}<p$, and hence the claim. Thus, we assume that $q+q_{n-1}\geq p$, and necessarily $h<a-a_n$. For $i=q-q_{n-1}$ and $j=q_{n-1}-q_{n-2}$, we have that $0\leq i,j\leq p-1$ and $i+j\leq p-1$. Indeed, it is clear that $j\geq0$ and $i+j\leq p-1$. Let us check that $i\geq0$. We have that $\widehat{q\frac{a}{p}}=\frac{p+h-a}{p}\leq\frac{p-1}{p}=\widehat{q_{n-1}\frac{a}{p}}$. By Proposition \ref{pro:contfrac} \ref{pro:contfrac2}, $q\geq q_{n-1}$. Now, $r_i=[(q-q_{n-1}+1)a-h]_p=1$ and $r_j=[(q_{n-1}-q_{n-2}+1)a-h]_p=a-a_n-h-1$. We get $$r_i+r_j=a-h-a_n<a-h,$$ and \eqref{eq:ineqremainders} is not satisfied.

Suppose now that $h=a-a_n$, and let us prove that \eqref{eq:ineqremainders} holds for all $0\leq i,j\leq p-1$ with $i+j\leq p-1$. Otherwise, suppose that there are $0\leq i,j\leq p-1$ with $i+j\leq p-1$ such that \eqref{eq:ineqremainders} does not hold; that is, $r_i+r_j<a_n$. Then, $r_i,r_j<a_n$. We shall characterize the integers $0\leq k\leq p-1$ such that $r_k<a_n$. By definition, $ka+a_n\equiv r_k\,(\mathrm{mod}\,p)$, so $ka\equiv-(a_n-r_k)\,(\mathrm{mod}\,p)$. Multiplying by $q_{n-1}$, we get $-k\equiv-(a_n-r_k)q_{n-1}\,(\mathrm{mod}\,p)$. Equivalently, $k\equiv(a_n-r_k)q_{n-1}\,(\mathrm{mod}\,p)$. Now, $(a_n-r_k)q_{n-1}<q_{n-1}a_n<p$, so $$k=(a_n-r_k)q_{n-1}.$$ In our situation, we find that $i=xq_{n-1}$ and $j=yq_{n-1}$ with $x=a_n-r_i$ and $y=a_n-r_j$. Note that $$a_n>r_i+r_j=2a_n-(x+y),$$ so $x+y>a_n$; that is, $x+y\geq a_n+1$. Therefore, $$i+j=(x+y)q_{n-1}\geq (a_n+1)q_{n-1}=a_nq_{n-1}+q_{n-1}\geq a_nq_{n-1}+q_{n-2}=p,$$ which contradicts that $i+j\leq p-1$.

Let us move on to the case in which $n$ is odd. In this case, we have that $[q_{n-1}a]_p=1$, so $q_{n-1}a\equiv1\,(\mathrm{mod}\,p)$ and $q_{n-2}a\equiv-a_n\,(\mathrm{mod}\,p)$.

First, we prove that if $0\leq h<a-\frac{a_n}{2}$, then \eqref{eq:ineqremainders} does not hold for some $0\leq i,j\leq p-1$ with $i+j\leq p-1$. Note that the previous inequality yields $a-h>\frac{a_n}{2}$. By definition of continued fraction, $a_n\geq2$, so $a-h>1$. If $q\leq\frac{p-1}{2}$, choose $i=j=q$. In that case, $i+j=2q\leq p-1$ and $0\leq i,j\leq p-1$. Now, $r_i=r_j=0$, so $r_i+r_j=0<a-h$. Next, suppose that $\frac{p-1}{2}<q<p-q_{n-1}$, and pick $i=q$, $j=q_{n-2}$. Again, $0\leq i,j\leq p-1$, and $i+j=q+q_{n-2}\leq q+q_{n-1}\leq p-1$. We already know that $r_i=0$. On the other hand, $r_j\equiv q_{n-2}a+a-h\equiv a-h-a_n\,(\mathrm{mod}\,p)$. We claim that $a-h>a_n$. Otherwise, assume that $a-h\leq a_n$. Since $qa\equiv-(a-h)\,(\mathrm{mod}\,p)$ and $q_{n-1}a\equiv1\,(\mathrm{mod}\,p)$, then $q\equiv-(a-h)q_{n-1}\,(\mathrm{mod}\,p)$, and $a-h\leq a_n$ implies that $(a-h)q_{n-1}\leq a_nq_{n-1}<p$, so $0<p-(a-h)q_{n-1}<p$. Therefore, $q=p-(a-h)q_{n-1}$. Now, since $a-h>\frac{a_n}{2}$, we have that $2(a-h)\geq a_n+1$, so $2(a-h)q_{n-1}\geq(a_n+1)q_{n-1}>p$; equivalently, $(a-h)q_{n-1}>\frac{p}{2}$. Hence, $q<p-\frac{p}{2}=\frac{p}{2}$, which contradicts that $q>\frac{p-1}{2}$. Hence, $r_j=a-h-a_n$, so $r_i+r_j=a-h-a_n<a-h$. Finally, assume that $q\geq p-q_{n-1}$. Let us choose $i=q+q_{n-1}-p$ and $j=p-2q_{n-1}$. Note that $i\geq0$ follows from our assumption and that $j\geq0$ because $a_n\geq2$. Moreover, $i+j=q-q_{n-1}\leq p-1$. Now, $$r_i=[(q+q_{n-1}-p)a+a-h]_p=[(q_{n-1}-p)a]_p=1.$$ On the other hand, $$r_j=[(p-2q_{n-1})a+a-h]_p=[a-h-2]_p=a-h-2.$$ Hence, $r_i+r_j=a-h-1<a-h$, proving that \eqref{eq:ineqremainders} fails.

Finally, in order to validate \ref{pro:characteqathetah2}, it remains to prove that if $a-\frac{a_n}{2}\leq h<a$, then \eqref{eq:ineqremainders} holds. Note that our assumption is equivalent to $1\leq a-h\leq\frac{a_n}{2}$. Suppose that there are $0\leq i,j\leq p-1$ with $i+j\leq p-1$ such that \eqref{eq:ineqremainders} fails, that is, $r_i+r_j<a-h$, whence $r_i,r_j<a-h$. We argue as in the analogous case for even $n$: if $0\leq k\leq p-1$ is such that $r_k<a-h$, then $ka+a-h\equiv r_k\,(\mathrm{mod}\,p)$, so $ka\equiv-(a-h-r_k)\,(\mathrm{mod}\,p)$. Multiplying by $q_{n-1}$, we obtain $k\equiv-(a-h-r_k)q_{n-1}\,(\mathrm{mod}\,p)$. Now, $$0<(a-h-r_k)q_{n-1}\leq(a-h)q_{n-1}<a_nq_{n-1}<p,$$ so $$k=p-(a-h-r_k)q_{n-1}.$$ On the other hand, multiplying by $q_{n-1}$ in $qa+a-h\equiv0\,(\mathrm{mod}\,p)$ yields $q+(a-h)q_{n-1}\equiv0\,(\mathrm{mod}\,p)$. In addition, $0<(a-h)q_{n-1}<p$, so $q=p-(a-h)q_{n-1}$. Plugging this in the expression of $k$ above, we get $k=q+r_kq_{n-1}$. Note that, since $(a-h)q_{n-1}<\frac{a_nq_{n-1}}{2}<\frac{p}{2}$, we have that $q>\frac{p}{2}$. Therefore, $k\geq q>\frac{p}{2}$. Applying this with $k=i$ and $k=j$ successively, we get that $i,j>\frac{p}{2}$. Hence, $i+j>p$, which is a contradiction. This finishes the proof.
\end{proof}

\begin{rmk}\normalfont For the case $h=0$, in \cite[Proposition 6.11]{gildegpgeneral} it is proved that $\mathfrak{A}_{L/K}=\mathfrak{A}_{\theta}$ if and only if $a\mid p-1$. This result can be recovered from \ref{pro:characteqathetah} as follows. It is straightforward to check that $a\mid p-1$ if and only if $n\leq 2$ (see for instance \cite[Proposition 6.12]{gildegpgeneral}). If $a=1$, then $a\mid p-1$ and from \ref{pro:characteqathetah1} we have that $\mathfrak{A}_{L/K}=\mathfrak{A}_{\theta}$. If $a>1$, from \ref{pro:characteqathetah2} we get:
\begin{itemize}
    \item If $n$ is even, $\mathfrak{A}_{L/K}=\mathfrak{A}_{\theta}$ if and only if $a=a_n$. Note that $a=p_n$. The numerators $\{p_i\}_{i=1}^n$ of the convergents satisfy the recurrence relation $p_i=a_ip_{i-1}+p_{i-2}$ for all $2\leq i\leq n$, so $a_n=p_n$ if and only if $n=2$, which is equivalent to $a\mid p-1$.
    \item If $n$ is odd, $\mathfrak{A}_{L/K}=\mathfrak{A}_{\theta}$ if and only if $a\leq\frac{a_n}{2}$. Again, since $a=p_n$, we have that $p_n\leq\frac{a_n}{2}$ if and only if $n=1$, which is equivalent to $a=1$, and this divides $p-1$.
\end{itemize} 
\end{rmk}

\subsection{The case $a+1\leq h\leq p-1$}

Under the assumption that $a+1\leq h\leq p-1$, we have that $a'=a+p$, so for every $0\leq k\leq p-1$, we have $$d(k;h)=\Big\lfloor\frac{a+kb-h}{p}\Big\rfloor+1=ka_0+\Big\lfloor\frac{(k+1)a-h}{p}\Big\rfloor+1.$$

\begin{pro}\label{pro:charactineqgreater} Suppose that $a+1\leq h\leq p-1$. Then, given $0\leq i,j\leq p-1$, the following statements hold:
\begin{enumerate}[label=\arabic*)]
    \item\label{pro:charactineqgreater1} If $i+j\leq p-1$, then $d(i;h)+d(j;h)\leq d(i+j;h)$ if and only if \begin{equation}\label{eq:ineqremainders2}
    r_i+r_j\geq p+a-h.
\end{equation}
    \item\label{pro:charactineqgreater2} If $i+j\geq p$ and $\ell<\frac{pe}{p-1}-1$, then $d(i;h)+d(j;h)\leq v_K(\lambda)+d(i+j+1-p;h)$.
\end{enumerate}
\end{pro}
\begin{proof}
Fix $0\leq i,j\leq p-1$.
\begin{enumerate}[label=\arabic*)]
    \item Suppose that $i+j\leq p-1$. Arguing as in the proof of Proposition \ref{pro:ineqremainders} \ref{pro:ineqremainders1}, we see that $d(i;h)+d(j;h)\leq d(i+j;h)$ if and only if \begin{equation}\label{eq:ineqfractpart2}
    \widehat{\frac{(i+1)a-h}{p}}+\widehat{\frac{(j+1)a-h}{p}}\geq\frac{a-h}{p}+\widehat{\frac{(i+j+1)a-h}{p}}+1.
\end{equation} This is clearly equivalent to $r_i+r_j\geq p+a-h+r_{i+j}$. Now, again we have $$r_i+r_j=\Big\lfloor\frac{r_i+r_j-(a-h)}{p}\Big\rfloor p+a-h+r_{i+j}.$$ Hence, our inequality is fulfilled if and only if $\Big\lfloor\frac{r_i+r_j-(a-h)}{p}\Big\rfloor\geq 1$, which is trivially equivalent to \eqref{eq:ineqremainders2}.
    \item Suppose that $i+j\geq p$ and $\ell<\frac{pe}{p-1}-1$. Since $a+1\leq h\leq p-1$, in particular, $a<p-1$. By Lemma \ref{lem:relprecramjump}, we obtain $\mathfrak{c}\geq p+a$, so $d(p-1)<v_K(\lambda)$ (using the expression of $\mathfrak{c}$ from Proposition \ref{pro:precscaffold}). Now, \begin{equation*}
    \begin{split}
        d(i;h)+d(j;h)&=(i+j)a_0+\Big\lfloor\frac{(i+1)a-h}{p}\Big\rfloor+\Big\lfloor\frac{(j+1)a-h}{p}\Big\rfloor+2
        \\&\leq(i+j)a_0+\Big\lfloor\frac{(i+j+2)a-h}{p}\Big\rfloor+2
        \\&=d(p-1)+\Big\lfloor\frac{a+(i+j+1-p)b-h}{p}\Big\rfloor+2
        \\&=d(p-1)+1+\Big\lfloor\frac{a'+(i+j+1-p)b-h}{p}\Big\rfloor
        \\&\leq v_K(\lambda)+d(i+j+1-p;h).
    \end{split}
\end{equation*}
\end{enumerate}
\end{proof}

\begin{pro} Suppose that $a+1\leq h\leq p-1$ and $\ell<\frac{pe}{p-1}-1$. Then, $\mathfrak{A}_{L/K}(h)=\mathfrak{A}_{\theta}(h)$ if and only if $a=1$ and $h>\frac{p+1}{2}$.
\end{pro}
\begin{proof}
The strategy is similar as in the proof of Proposition \ref{pro:characteqathetah}. Namely, by Corollary \ref{coro:twoimpl} and Proposition \ref{pro:charactineqgreater}, $\mathfrak{A}_{L/K}(h)=\mathfrak{A}_{\theta}(h)$ if and only if \eqref{eq:ineqremainders2} holds. Once again, we call $q\coloneqq[a^{-1}h-1]_p$. 

First, assume that $a=1$ (recall that this is equivalent to $n=1$). Suppose that $h\leq\frac{p+1}{2}$, so that $h-a\leq\frac{p-1}{2}$. Choose $i=j=h-a$. Then $i,j>0$ and $i+j=2(h-a)\leq p-1$. Now, $r_i=r_j=[h-a+a-h]_p=0$, so $r_i+r_j=0<p+a-h$, and \eqref{eq:ineqremainders2} does not hold. Otherwise, we assume $h>\frac{p+1}{2}$, so that $h-a>\frac{p-1}{2}$. Given $0\leq k\leq p-1$, if $k<h-a$ then $$r_k=[k+a-h]_p=p+k+a-h\geq p+a-h.$$ Let $0\leq i,j\leq p-1$. We apply the previous argument with $k=i,j$, to obtain that if $r_i+r_j<p+a-h$, then $r_i,r_j<p+a-h$. It follows that $i,j\geq h-a$, and hence $i+j\geq2(h-a)>p-1$. Equivalently, $i+j\leq p-1$ implies that $r_i+r_j\geq p+a-h$, so \eqref{eq:ineqremainders2} holds. 

Next, we assume that $a>1$. If $q=p-1$, then $h-a\equiv p-a\,(\mathrm{mod}\,p)$ with $1\leq h-a,p-a<p$, so $h=p$, which is a contradiction. Hence, $q\leq p-2$. Suppose that $h\leq2a$, and choose $i=q$, $j=1$. Then $0\leq i,j\leq p-1$ and $i+j=q+1\leq p-1$. Now, $r_i=0$ and $r_j=[2a-h]_p=2a-h$, as $h>a$ implies $2a-h<2a-a=a$. We see that $r_i+r_j<p+a-h$ and \eqref{eq:ineqremainders2} is not satisfied. 

Let us assume now that $h>2a$. Define $m\coloneqq\Big\lceil\frac{h-a}{a}\Big\rceil$. It is straightforward to check that $(m-1)a<h-a\leq ma$. Set $\delta\coloneqq ma+a-h$. Then, $0\leq\delta<ma-(m-1)a=a$, that is, $0\leq\delta\leq a-1$. In particular, $r_m=[\delta]_p=\delta$. Also, $ma=\delta+h-a\leq a-1+h-a=h-1\leq p-2$. Since $a\geq 2$, we deduce that $2m<p-1$. 

Suppose that $q+m\leq p-1$, and let us choose $i=q$ and $j=m$. Then $r_i+r_j=r_m=\delta\leq a-1$. On the other hand, $p+a-h\geq p+a-(p-1)=a+1$. We see that $r_i+r_j<p+a-h$, so \eqref{eq:ineqremainders2} does not hold. 

Finally, suppose that $q+m\geq p$. We have $$(p-q)a\equiv a-h\equiv p+a-h\,(\mathrm{mod}\,p).$$ Our assumption together with $q\leq p-2$ yield $2\leq p-q\leq m$. Hence, $1\leq(p-q)a,p+a-h\leq p-1$. We obtain that $(p-q)a=p+a-h$. Let $i=j=m$. Then $0\leq i,j\leq p-1$ with $i+j=2m<p-1$. Now, $$r_i+r_j=2\delta<2a\leq(p-q)a=p+a-h,$$ disproving \eqref{eq:ineqremainders2}. This concludes the proof.
\end{proof}

\section{Broadening the picture}\label{sec:exceptionalcase}

In this part, we shall deal with the only uncovered case by Theorem \ref{thm:freenessideals}: $a+1\leq h\leq p-1$ and $\ell\geq\frac{pe}{p-1}-1$. This is the so called exceptional case. We shall prove Theorem \ref{thm:assocordtypical} on a description of $\mathfrak{A}_{L/K}(h)$ valid for all typical degree $p$ extensions; in particular for the ones at the exceptional case. Afterwards, we will specialize to this case and prove Theorem \ref{thm:freenessexceptionalcase}.

If $0\leq h\leq a$ or $\ell<\frac{pe}{p-1}-1$, the parameters $w(i;h)$ provided a big deal of the information on both the associated order $\mathfrak{A}_{L/K}(h)$ and the freeness of $\mathfrak{p}_L^h$. In the general situation, we shall present new parameters $\widetilde{w}(i;h)$ that will improve the role of the $w(i;h)$, in the following sense. In the already known cases, we will have $\widetilde{w}(i;h)=w(i;h)$. Moreover, in all those cases that the $w(i;h)$ are not enough so as to describe the associated order and the freeness, we will show that the $\widetilde{w}(i;h)$ do so.

Fix $0\leq i\leq p-1$ and recall that $w(i;h)=\underset{0\leq j\leq p-1-i}{\mathrm{min}}W(i,j;h)$, where $W(i,j;h)=d(i+j;h)-d(j;h)$. Note that $W(i,0;h)=d(i;h)$, which is why $w(i;h)\leq d(i;h)$. Let $\theta=\pi_L^{a'}$. By Corollary \ref{coro:twoimpl}, $\mathfrak{A}_{\theta}(h)$ being a ring is equivalent to the implications \ref{coro:twoimpl1} and \ref{coro:twoimpl2}:
\begin{itemize}
    \item If $i+j\leq p-1$, then $d(i;h)+d(j;h)\leq d(i+j;h)$.
    \item If $i+j\geq p$, then $d(i;h)+d(j;h)\leq v_K(\lambda)+d(i+j+1-p;h)$.
\end{itemize}

Under the assumption that $0\leq h\leq a$ or $\ell<\frac{pe}{p-1}$, we have seen in Propositions \ref{pro:ineqremainders} \ref{pro:ineqremainders2} and \ref{pro:charactineqgreater} \ref{pro:charactineqgreater2} that the second one holds, so the freeness is characterized by the first one, which is equivalent to $d(i;h)\leq W(i,j;h)$. If now we move on to the general case, the second implication has to be considered as well. Our course of action is extending the minimum defining the parameters $w(i;h)$ to all $0\leq j\leq p-1$.

Given $0\leq i,j\leq p-1$, let us denote $$\varepsilon(i,j)=\begin{cases}
    0 & \hbox{if }i+j\leq p-1, \\
    1 & \hbox{if }i+j\geq p,
\end{cases}$$ $$i\oplus j=\begin{cases}
    i+j & \hbox{if }i+j\leq p-1, \\
    i+j+1-p & \hbox{if }i+j\geq p.
\end{cases}$$ Let us call $$W(i,j;h)\coloneqq\varepsilon(i,j)v_K(\lambda)+d(i\oplus j;h)-d(j;h)$$ for all such $i$ and $j$. If $0\leq j\leq p-1-i$, we have that $W(i,j;h)=d(i+j;h)-d(j;h)$, so this definition for $W(i,j;h)$ is compatible with the already existing one. Now, we define $$\widetilde{w}(i;h)=\mathrm{min}\{W(i,j;h)\,\mid\,0\leq j\leq p-1\},\quad 0\leq i\leq p-1.$$ When $h=0$, we shall call $\widetilde{w}(i)=\widetilde{w}(i;0)$ and $W(i,j)=W(i,j;0)$.

\begin{pro}\label{pro:propertparam} The following properties are satisfied.
\begin{enumerate}[label=\arabic*)]
    \item\label{pro:propertparam1} $\widetilde{w}(i;h)\leq w(i;h)\leq d(i;h)$ for all $0\leq i\leq p-1$.
    \item\label{pro:propertparam2} $\mathfrak{A}_{\theta}(h)$ is a ring if and only if $\widetilde{w}(i;h)=d(i;h)$ for all $0\leq i\leq p-1$, where $\theta=\pi_L^{a'}$.
    \item\label{pro:propertparam3} If $0\leq h\leq a$ or $\ell<\frac{pe}{p-1}-1$, then $\widetilde{w}(i;h)=w(i;h)$ for all $0\leq i\leq p-1$.
\end{enumerate}
\end{pro}
\begin{proof}
\begin{enumerate}
    \item We already know that $w(i;h)\leq d(i;h)$ for all $0\leq i\leq p-1$. The first inequality follows directly from the definition: Given $0\leq i\leq p-1$, if $0\leq j\leq p-1-i$ (so $i+j\leq p-1$), then $\varepsilon(i,j)=0$ and $i\oplus j=i+j$, and hence $W(i,j;h)=d(i+j;h)-d(j;h)$. Since $w(i;h)=\underset{0\leq j\leq p-1-i}{\mathrm{min}}(d(i+j;h)-d(j;h))$, either the minimum defining $\widetilde{w}(i;h)$ is reached at $0\leq j\leq p-1-i$, in which case $\widetilde{w}(i;h)=w(i;h)$, or else it is reached at $p-i\leq j\leq p-1$, in which case $\widetilde{w}(i;h)\leq w(i;h)$.
    \item By Corollary \ref{coro:twoimpl}, $\mathfrak{A}_{\theta}(h)$ is a ring if and only if $d(i;h)+d(j;h)\leq\varepsilon(i,j)v_K(\lambda)+d(i\oplus j;h)$ for all $0\leq i,j\leq p-1$. Equivalently, $d(i;h)\leq W(i,j;h)$ for all $0\leq i,j\leq p-1$, which holds if and only if $d(i;h)\leq\widetilde{w}(i;h)$ for all $0\leq i\leq p-1$. On the other hand, we know from \ref{pro:propertparam1} that $\widetilde{w}(i;h)\leq d(i;h)$ for all $0\leq i\leq p-1$. Thus, $\mathfrak{A}_{\theta}(h)$ is a ring if and only if $\widetilde{w}(i;h)=d(i;h)$ for all $0\leq i\leq p-1$, as we wanted.
    \item Fix $0\leq i\leq p-1$ and $p-i\leq j\leq p-1$. By Propositions \ref{pro:ineqremainders} \ref{pro:ineqremainders2} and \ref{pro:charactineqgreater} \ref{pro:charactineqgreater2}, under the stated hypotheses we have that $d(i;h)+d(j;h)\leq v_K(\lambda)+d(i+j+1-p;h)$. This is the same as $d(i;h)\leq W(i,j;h)$. But recall that $d(i;h)=W(i,0;h)$, so the minimum has to be reached at some $0\leq j'\leq p-1-i$. This proves that $\widetilde{w}(i;h)=w(i;h)$.
\end{enumerate}
\end{proof}

We consider the problem of describing $\mathfrak{A}_{L/K}(h)$.

\begin{lem}\label{lem:testassocord} Let $L/K$ be a typical degree $p$ extension of $p$-adic fields with the usual notation. Call $\theta=\pi_L^{a'}$. An element $\gamma=\sum_{i=0}^{p-1}\gamma_i\Psi^i\in H$ with $\gamma_i\in K$ for all $0\leq i\leq p-1$ belongs to $\mathfrak{A}_{L/K}(h)$ if and only if each summand $\gamma_i\Psi^i$ does.
\end{lem}
\begin{proof}
We use repeatedly Corollary \ref{coro:quotideal}, by which an element $\rho\in H$ belongs to $\mathfrak{A}_{L/K}$ if and only if $\rho\mathfrak{A}_{\theta}(h)\subseteq\mathfrak{A}_{\theta}(h)$. By Proposition \ref{pro:basisathetah}, this happens if and only if $\rho f_j\in\mathfrak{A}_{\theta}(h)$ for all $0\leq j\leq p-1$, and this is clearly equivalent to $\rho f_j\cdot\theta\in\mathfrak{p}_L^h$ for all $0\leq j\leq p-1$.

The right to left implication is trivial. Conversely, suppose that $\gamma\in\mathfrak{A}_{L/K}(h)$. Then, $\gamma f_j\cdot\theta\in\mathfrak{p}_L^h$ for all $0\leq j\leq p-1$. In particular, for $j=0$, we obtain $\gamma\cdot\theta\in\mathfrak{p}_L^h$ (since $f_0=\mathrm{Id}$). Now, $$\gamma\cdot\theta=\sum_{i=0}^{p-1}\gamma_i\Psi^i\cdot\theta$$ and for each $0\leq i\leq p-1$ we have $$v_L(\gamma_i\Psi^i\cdot\theta)=a'+ib\equiv(i+1)a\,(\mathrm{mod}\,p).$$ Since these valuations are distinct mod $p$, we obtain that $\gamma_i\Psi^i\cdot\theta\in\mathfrak{p}_L^h$ for all $0\leq i\leq p-1$. In particular, for $i=0$, we have that $\gamma_0\theta\in\mathfrak{p}_L^h$. Taking valuations, $h\leq v_L(\gamma_0\theta)=pv_K(\gamma_0)+a'$. Since $h\leq a'\leq h+p-1$, $$v_K(\gamma_0)\geq\Big\lceil\frac{h-a'}{p}\Big\rceil=0,$$ whence $\gamma_0\in\mathcal{O}_K$, and hence $\gamma_0\mathrm{Id}\in\mathfrak{A}_{L/K}$. It remains to check it for $1\leq i\leq p-1$. From the case $i=0$ we have that $$\gamma-\gamma_0\mathrm{Id}=\sum_{i=1}^{p-1}\gamma_i\Psi^i\in\mathfrak{A}_{L/K}(h).$$ Consequently, for every $0\leq j\leq p-1$, we have that $$\sum_{i=1}^{p-1}\gamma_i\Psi^if_j\cdot\theta\in\mathfrak{p}_L^h.$$ Since $$\Psi^if_j=\pi_K^{-d(j;h)}\lambda^{\varepsilon(i,j)}\Psi^{i\oplus j},$$ for each $1\leq i\leq p-1$ we have $$v_L(\gamma_i\Psi^if_j\cdot\theta)\equiv a'+(i\oplus j)b\equiv(i+j+1)a\,(\mathrm{mod}\,p).$$ These $L$-valuations are pairwise distinct, proving that $$\gamma_i\Psi^if_j\cdot\theta\in\mathfrak{p}_L^h$$ for every $1\leq i\leq p-1$ and every $0\leq j\leq p-1$. Hence, $\gamma_i\Psi^i\in\mathfrak{A}_{L/K}(h)$ for all $1\leq i\leq p-1$, finishing the proof.
\end{proof}

For each $0\leq i\leq p-1$, let us define $$\widetilde{e}_i=\pi_K^{-\widetilde{w}(i;h)}\Psi^i.$$ We restate and prove Theorem \ref{thm:assocordtypical}.

\begin{pro}\label{pro:assocordalmostmaxram} Let $L/K$ be a typical degree $p$ extension of $p$-adic fields with the usual notation. The associated order $\mathfrak{A}_{L/K}(h)$ admits an $\mathcal{O}_K$-basis given by $\{\widetilde{e}_i\}_{i=0}^{p-1}$.
\end{pro}
\begin{proof}
Pick $\gamma=\sum_{i=0}^{p-1}\gamma_i\Psi^i\in H$. By Lemma \ref{lem:testassocord}, $\gamma\in\mathfrak{A}_H(h)$ if and only if $\gamma_i\Psi^i\in\mathfrak{A}_{L/K}(h)$ for all $0\leq i\leq p-1$.

Fix $0\leq i\leq p-1$. We already know that $\gamma_i\Psi^i\in\mathfrak{A}_{L/K}(h)$ if and only if $\gamma_i\Psi^if_j\in\mathfrak{A}_{\theta}(h)$ for all $0\leq j\leq p-1$. We also fix some such a $j$. Now, \begin{equation*}
    \begin{split}
        \gamma_i\Psi^if_j&=\gamma_i\pi_K^{-d(j;h)}\Psi^i\Psi^j
        \\&=\gamma_i\pi_K^{-d(j;h)}\lambda^{\varepsilon(i,j)}\Psi^{i\oplus j}.
    \end{split}
\end{equation*} This belongs to $\mathfrak{A}_{\theta}(h)$ if and only if $$p(v_K(\gamma_i)-d(j;h)+\varepsilon(i,j)v_K(\lambda))+a'+(i\oplus j)b\geq h.$$ Subtracting $h$ and dividing by $p$ gives $$v_K(\gamma_i)-d(j;h)+\varepsilon(i,j)v_K(\lambda)+d(i\oplus j;h)\geq0.$$ Finally, reorganizing suitably the terms yields $v_K(\gamma_i)\geq-W(i,j;h)$. Since $j$ is arbitrary, this is equivalent to $v_K(\gamma_i)\geq-\widetilde{w}(i;h)$. Therefore, $\gamma\in\mathfrak{A}_{L/K}(h)$ if and only if $\gamma_i\in\pi_K^{-\widetilde{w}(i;h)}\mathcal{O}_K$ for all $0\leq i\leq p-1$, proving the statement.
\end{proof}

The proof of Proposition \ref{pro:assocordalmostmaxram} contains as a particular case a proof that the elements $e_i=\pi_K^{-w(i;h)}\Psi^i$ form an $\mathcal{O}_K$-basis of $\mathfrak{A}_{L/K}(h)$ if $0\leq h\leq a$ or $\ell<\frac{pe}{p-1}-1$ (which was already established in Corollary \ref{coro:criteriamidtransl} \ref{coro:criteriamidtransl1}). Indeed, in that case, if $\gamma\in\mathfrak{A}_{L/K}(h)$, then $\gamma_i\Psi^if_j\in\mathfrak{A}_{\theta}(h)$ always holds when $i+j\geq p$, because $\gamma_i\Psi^if_j\in\mathfrak{A}_{\theta}(h)$ is equivalent to $d(i;h)+d(j;h)\leq v_K(\lambda)+d(i+j+1-p;h)$, which holds without further restrictions by Propositions \ref{pro:ineqremainders} \ref{pro:ineqremainders2} and \ref{pro:charactineqgreater} \ref{pro:charactineqgreater2}. If $i+j\leq p-1$, the condition that $\gamma_i\Psi^if_j\in\mathfrak{A}_{\theta}(h)$ ends up reducing to $v_K(\gamma_i)\geq-w(i;h)$. In fact, this argument appears in \cite[Theorem 7.2]{delcorsoferrilombardo} for the ring of integers in the Galois case.

\subsection{A characterization for freeness}

In this part, we assume that $a+1\leq h\leq p-1$ and $\ell\geq\frac{pe}{p-1}-1$. Under these conditions, we know that $a'=a+p$ and $d(i;h)=\Big\lfloor\frac{a+ib-h}{p}\Big\rfloor+1$.

Let us characterize when $\mathfrak{p}_L^h$ is $\mathfrak{A}_{L/K}(h)$-free. Let us call again $A_h=\mathfrak{A}_{L/K}(h)$ and $I_h=\mathfrak{A}_{\theta}(h)$, and set $$P_1\coloneqq\lambda^{-1}\Psi^{p-1},\quad P_0\coloneqq1-P_1.$$ Note that $$P_1^2=\lambda^{-2}\Psi^{2p-2}=\lambda^{-1}\Psi^{p-1}=P_1,$$ so $P_1$ is an idempotent of $H$, and then, so is $P_0$. Also, $P_1\Psi^j=\Psi^j$ and $P_0\Psi^j=0$ for all $1\leq j\leq p-1$. Hence $P_1I_h\subseteq I_h$ and $P_0I_h\subseteq I_h$, so $P_1,P_0\in A_h$ by Corollary \ref{coro:quotideal}.

We shall determine decompositions of $A_h$ and $I_h$ in terms of these idempotents. Let $H^{(1)}=P_1H$. Note that this is a ring with identity $P_1$. Now, let us denote $A_h^{(1)}=P_1A_h$ and $I_h^{(1)}=P_1I_h$. Since $P_1$ is a central idempotent, $A_h^{(1)}$ is a subring of $H^{(1)}$. On the other hand, since $P_0\Psi^j=0$ for $1\leq j\leq p-1$, $P_0A_h=\mathcal{O}_KP_0$. Thus $P_1+P_0=1$ yields $$A_h=\mathcal{O}_KP_0\oplus A_h^{(1)}.$$ Similarly, $$I_h=\mathcal{O}_KP_0\oplus I_h^{(1)}.$$ Let us see that $I_h^{(1)}$ is an $A_h^{(1)}$-fractional ideal. Since $I_h$ is an $A_h$-fractional ideal, there is $\alpha\in A_h$ such that $\alpha I_h\subseteq A_h$. Multiplying by $P_1$ and taking into account that $P_1=P_1^2$, we get $(P_1\alpha) I_h^{(1)}\subseteq A_h^{(1)}$, which proves the claim. We derive the following simplification of our problem.

\begin{lem} $I_h$ is $A_h$-principal if and only if $I_h^{(1)}$ is $A_h^{(1)}$-principal.
\end{lem}
\begin{proof}
The implication from left to right is immediate by multiplication with $P_1$. Conversely, suppose that $I_h^{(1)}$ is $A_h^{(1)}$-principal, and let $\gamma_1\in I_h^{(1)}$ be such that $I_h^{(1)}=A_h^{(1)}\gamma_1$. Let $\gamma\coloneqq P_0+\gamma_1$, which belongs to $I_h$ because $I_h=P_0\mathcal{O}_K\oplus I_h^{(1)}$. We have \begin{equation*}
    \begin{split}
        A_h\gamma&=(P_0\mathcal{O}_K\oplus A_h^{(1)})(P_0+\gamma_1)
        \\&=P_0\mathcal{O}_K\oplus A_h^{(1)}\gamma_1
        \\&=P_0\mathcal{O}_K\oplus I_h^{(1)}=I_h,
    \end{split}
\end{equation*} whence $I_h$ is $A_h$-principal.
\end{proof}

Thus, to study the principality (and hence the freeness), we can ignore the $P_0$-component. But at the level of the already known bases for $I_h$ and $A_h$, the $P_0$-component is equivalent to the identity component. Explicitly:

\begin{pro} An $\mathcal{O}_K$-basis of $A_h$ (resp. $I_h$) consists in the element $P_0$ together with the elements $\widetilde{e}_i$ (resp. $f_i$), where $1\leq i\leq p-1$.
\end{pro}
\begin{proof}
It is enough to prove that $P_1$ and $\widetilde{e}_{p-1}$ coincide up to multiplication by a unit of $\mathcal{O}_K$, whence it will follow immediately that $1$ can be replaced by $P_0=1-P_1$ in the corresponding basis. We have that $$P_1=\lambda^{-1}\Psi^{p-1}=\lambda^{-1}\pi_K^{\widetilde{w}(p-1;h)}\widetilde{e}_{p-1}$$ with $v_K(\lambda^{-1}\pi_K^{\widetilde{w}(p-1;h)})=\widetilde{w}(p-1;h)-v_K(\lambda)$, and this vanishes. Indeed, we find $\widetilde{w}(p-1;h)$ as the minimum on $j$ of $W(i,j;h)$. For $j=0$, we have $$W(p-1,0;h)=\varepsilon(p-1,0)v_K(\lambda)+d(p-1;h)-d(0;h)=d(p-1;h),$$ and $d(p-1;h)=d(p-1)$ follows as in the proof of Lemma \ref{lem:diffequal1}. Since $\ell\geq\frac{pe}{p-1}-1$, we have that $d(p-1)=v_K(\lambda)$. On the other hand, for $j>0$, we get $$W(p-1,j;h)=\varepsilon(p-1,j)v_K(\lambda)+d(j;h)-d(j;h)=v_K(\lambda).$$
\end{proof}

Equivalently, we have that $A_h=\mathcal{O}_KP_0\oplus\Big(\bigoplus_{i=1}^{p-1}\mathcal{O}_K\widetilde{e}_i\Big)$ and $I_h=\mathcal{O}_KP_0\oplus\Big(\bigoplus_{i=1}^{p-1}\mathcal{O}_Kf_i\Big)$. It follows that $$A_h^{(1)}=\bigoplus_{i=1}^{p-1}\mathcal{O}_K\widetilde{e}_i,\quad I_h^{(1)}=\bigoplus_{i=1}^{p-1}\mathcal{O}_Kf_i.$$

\begin{thm} Suppose that $a+1\leq h\leq p-1$ and $\ell\geq\frac{pe}{p-1}-1$. Then $\mathfrak{p}_L^h$ is $A_h$-free if and only if there is some $0\leq s\leq p-1$ such that $\widetilde{w}(i;h)=W(i,s;h)$ for all $0\leq i\leq p-1$. In other words, the same value $s$ achieves the minimum defining $\widetilde{w}(i;h)$ for all $0\leq i\leq p-1$.
\end{thm}
\begin{proof}
In order to ease the notation, we shall call $A=A_h^{(1)}$ and $I=I_h^{(1)}$ throughout this proof. We already know that $\mathfrak{p}_L^h$ is $A_h$-free if and only if $I$ is $A$-principal. In addition, if we call $A_i=\mathcal{O}_K\widetilde{e}_i$ and $I_i=\mathcal{O}_Kf_i$ for all $1\leq i\leq p-1$, then $$A=\bigoplus_{i=1}^{p-1}A_i,\quad I=\bigoplus_{i=1}^{p-1}I_i.$$ Analogously, $$H^{(1)}=\bigoplus_{i=1}^{p-1}K\Psi^i.$$

Suppose that for some $0\leq s\leq p-1$ we have $$\widetilde{w}(i;h)=W(i,s;h)=\varepsilon(i,s)v_K(\lambda)+d(i\oplus s;h)-d(s;h)$$ for all $0\leq i\leq p-1$. For every such an $i$, we have that $$\widetilde{e}_if_s=\pi_K^{-\widetilde{w}(i;h)-d(s;h)}\lambda^{\varepsilon(i,s)}\Psi^{i\oplus s}.$$ Now, $$v_K(\pi_K^{-\widetilde{w}(i;h)-d(s;h)}\lambda^{\varepsilon(i,s)})=-\widetilde{w}(i;h)-d(s;h)+\varepsilon(i,s)v_K(\lambda)=-d(i\oplus s;h),$$ the last equality due to our hypothesis. Hence, $\widetilde{e}_if_s$ coincides with $f_{i\oplus s}$ up to multiplication by a unit. Since $W(0,j;h)=0$ for all $0\leq j\leq p-1$ and $i\mapsto i\oplus s$ is a permutation of $\{1,\dots,p-1\}$, we see that $I=Af_s$, so $I$ is $A$-principal with generator $f_s$.

Conversely, assume that $I$ is $A$-principal, say $I=A\alpha$ for some $\alpha\in I$. Tensoring by $K$ gives $H^{(1)}=H^{(1)}\alpha$. In particular, $P_1\in H^{(1)}\alpha$, so there is $\alpha'\in H^{(1)}$ such that $P_1=\alpha\alpha'$. Since $P_1$ is the identity of $H^{(1)}$, $\alpha$ is a unit of $H^{(1)}$. Then, $I$ is invertible as a fractional ideal. Indeed, let us define $$J\coloneqq\{\gamma\in H^{(1)}\,\mid\,\gamma I\subseteq A\}.$$ Since $I=\alpha A$, necessarily $J=\alpha^{-1}A$, whence $IJ=A$.

Let us prove that $$J=\bigoplus_{i=1}^{p-1}J_i,\quad J_i=J\cap K\Psi^i.$$ Equivalently, given $\gamma\in J$, if we write $\gamma=\sum_{i=1}^{p-1}\gamma_i$ with $\gamma_i\in K\Psi^i$, then $\gamma_i\in J$ for all $1\leq i\leq p-1$. Fix some $1\leq j\leq p-1$ and take $\beta_j\in I_j$. Since $\gamma\in J$, by definition, $\gamma\beta_j\in A$. But $\gamma\beta_j=\sum_{i=1}^{p-1}\gamma_i\beta_j$, with $\gamma_i\beta_j\in K\Psi^{i\oplus j}$ for all $1\leq i\leq p-1$. We deduce that $\gamma_i\beta_j\in A_{i\oplus j}\subseteq A$ for all $1\leq i\leq p-1$. Since $j$ is arbitrary, $\gamma_i I\subseteq A$, whence $\gamma_i\in J$ for all $1\leq i\leq p-1$, as we wanted.

For each $1\leq j\leq p-1$, denote $$j^{\vee}=\begin{cases}
    p-1-j & \hbox{if }1\leq j<p-1, \\
    p-1 & \hbox{if }j=p-1.
\end{cases}$$ Note that $\{1,\dots,p-1\}$ endowed with the operation $(i,j)\mapsto i\oplus j$ is a group with identity $p-1$, whose inverse operation corresponds to $j\mapsto j^{\vee}$. Now, $$A_{p-1}=(IJ)_{p-1}=\sum_{j=1}^{p-1}I_jJ_{j^{\vee}},$$ but also $$A_{p-1}=\mathcal{O}_K\widetilde{e}_{p-1}=\mathcal{O}_KP_1,$$ as $P_1$ and $\widetilde{e}_{p-1}$ coincide up to multiplication by a unit of $\mathcal{O}_K$. For each $1\leq j\leq p-1$, let $n_j\in\mathbb{Z}_{\geq0}$ be such that $$I_jJ_{j^{\vee}}=\pi_K^{n_j}A_{p-1}.$$ Then, $$A_{p-1}=\sum_{j=1}^{p-1}\pi_K^{n_j}A_{p-1}=\pi_K^{\mathrm{min}\{n_j\,\mid\,1\leq j\leq p-1\}}A_{p-1}.$$ We deduce that $\mathrm{min}\{n_j\,\mid\,1\leq j\leq p-1\}=0$; that is, there is some $1\leq s\leq p-1$ such that $n_s=0$. Then $I_sJ_{s^{\vee}}=A_{p-1}=\mathcal{O}_KP_1$.

We claim that $I=Af_s$. It is trivial that $Af_s\subseteq I$. Conversely, let us take $\gamma\in I$. Let $\gamma_{s^{\vee}}$ be a generator of $J_{s^{\vee}}$ as an $\mathcal{O}_K$-module. Then $\gamma_{s^{\vee}}\gamma\in A$. On the other hand, $f_s\gamma_{s^{\vee}}$ is a generator of $I_sJ_{s^{\vee}}=\mathcal{O}_KP_1$, so there is $u\in\mathcal{O}_K^{\times}$ such that $f_s\gamma_{s^{\vee}}=uP_1$. Hence, $\gamma=(uP_1)^{-1}f_s\gamma_{s^{\vee}}\gamma=(u^{-1}\gamma_{s^{\vee}}\gamma)f_s$, where we have used that $P_1$ is the identity of $H^{(1)}$. Now, $u^{-1}\gamma_{s^{\vee}}\gamma\in JI=A$, so $\gamma\in Af_s$, and the desired equality follows.

Finally, we prove that $W(i,s;h)=\widetilde{w}(i;h)$ for all $0\leq i\leq p-1$. For $i=0$, it is immediate, as both sides of the desired equality vanish. Given $1\leq i\leq p-1$, we have that $A_if_s\subseteq A_iI_s\subseteq A_{i\oplus s}$. But precisely because that $I=Af_s$ and $i\mapsto i\oplus s$ is a permutation of $\{1,\dots,p-1\}$, we obtain that $A_if_s=I_{i\oplus s}$. Now, \begin{equation*}
    \begin{split}
        A_if_s&=\langle\pi_K^{-\widetilde{w}(i;h)}\Psi^i\rangle_{\mathcal{O}_K}(\pi_K^{-d(s;h)}\Psi^s)
        \\&=\langle\lambda^{\varepsilon(i,s)}\pi_K^{-\widetilde{w}(i;h)-d(s;h)}\Psi^{i\oplus s}\rangle_{\mathcal{O}_K}
        \\&=\langle\pi_K^{v_K(\lambda)\varepsilon(i,s)-\widetilde{w}(i;h)-d(s;h)}\Psi^{i\oplus s}\rangle_{\mathcal{O}_K}
    \end{split}
\end{equation*} On the other hand, $$I_{i\oplus s}=\langle\pi_K^{-d(i\oplus s;h)}\Psi^{i\oplus s}\rangle_{\mathcal{O}_K}.$$ Thus, the equality $A_if_s=I_{i\oplus s}$ yields $$v_K(\lambda)\varepsilon(i,s)-\widetilde{w}(i;h)-d(s;h)=-d(i\oplus s;h)$$ for all $0\leq i\leq p-1$. Equivalently, $\widetilde{w}(i;h)=W(i,s;h)$ for all $0\leq i\leq p-1$, finishing the proof.
\end{proof}

\section*{AI declaration}

During the preparation of this paper, the author used OpenAI’s ChatGPT (models GPT5.6 and GPT6). Before using these tools, the author had structured the project, defined its objectives, and established a strategy for pursuing them. ChatGPT assisted in developing results and proofs, particularly in Sections \ref{sec:equaldiwi} and \ref{sec:exceptionalcase}, and to a lesser extent in Section \ref{sec:equalassocorders}. It was also used to help revise the mathematical exposition and language throughout the paper. The author declares that has written every sentence of the manuscript and takes full responsability for the results, proofs, references, and final text.

\printbibliography

\end{document}